\documentclass[letterpaper, 10 pt, conference]{ieeeconf}
\IEEEoverridecommandlockouts 
\usepackage[table]{xcolor}
\usepackage{hyperref}
\usepackage{graphicx}
\usepackage{amsmath,amssymb}
\usepackage{cite}
\usepackage{float}
\usepackage{array}
\usepackage[caption=false,font=footnotesize]{subfig}
\usepackage{todonotes}
\usepackage{caption}
\usepackage{graphicx}

\newtheorem{Theorem}{Theorem}
\newtheorem{Definition}{Definition}
\newtheorem{Proposition}{Proposition}
\newtheorem{Lemma}{Lemma}

\hypersetup{
    colorlinks=true,
    citecolor=red,
    linkcolor=red,
    filecolor=green,      
    urlcolor=black,
}

\begin{document}

\renewcommand{\vec}{\boldsymbol}

\graphicspath{{img/},{fig/}}

%\markboth{aaa}{bbb}

\title{\LARGE \bf
Dimensionless Policies Based on the Buckingham $\pi$ Theorem: \\
Is This a Good Way to Generalize Numerical Results?
}

%\author{Alexandre~Girard,~\IEEEmembership{Member,~IEEE,}
%        and~H.~Harry~Asada,~\IEEEmembership{Member,~IEEE}% <-this % stops a space
\author{Alexandre Girard$^{1}$% <-this % stops a space
\thanks{$^{1}$Alexandre Girard is with the Department of Mechanical Engineering, Universite de Sherbrooke, Qc, Canada {\tt\small  alex.girard@usherbrooke.ca }}% <-this % stops a space
}%

% make the title area
\maketitle
\thispagestyle{empty}
\pagestyle{empty}

\begin{abstract}
    The answer to the question posed in the title is yes if the context (the list of variables defining the motion control problem) is dimensionally similar. This article explores the use of the Buckingham $\pi$ theorem as a tool to encode the control policies of physical systems into a more generic form of knowledge that can be reused in various situations. This approach can be interpreted as enforcing invariance to the scaling of the fundamental units in an algorithm learning a control policy. First, we show, by restating the solution to a motion control problem using dimensionless variables, that (1) the policy mapping involves a reduced number of parameters and (2) control policies generated numerically for a specific system can be transferred exactly to a subset of dimensionally similar systems by scaling the input and output variables appropriately. Those two generic theoretical results are then demonstrated, with numerically generated optimal controllers, for the classic motion control problem of swinging up a torque-limited inverted pendulum and positioning a vehicle in slippery conditions. We also discuss the concept of regime, a region in the space of context variables, that can help to relax the similarity condition. Furthermore, we discuss how applying dimensional scaling of the input and output of a context-specific black-box policy is equivalent to substituting new system parameters in an analytical equation under some conditions, using a linear quadratic regulator (LQR) and a computed torque controller as examples. It remains to be seen how practical this approach can be to generalize policies for more complex high-dimensional problems, but the early results show that it is a promising transfer learning tool for numerical approaches like dynamic programming and reinforcement learning.
\end{abstract}

%%%%%%%%%%%%%%%%%%%%%%%%%%%%%%%%%%%%%%%%%%%%
\section{Introduction}

% \paragraph{Many numerical algorithms = black box mapping}
% -Trajectory optimization
% -Reinforcement learning
% -etc.

% \paragraph{With numerical results, unlike analytical solutions, system and problem parameters are not explicitly in the solution }

% \paragraph{This makes the results "context specific" and makes its harder to generalize }

To solve challenging motion control problems in robotics (locomotion, manipulation, vehicle control, etc.), many approaches now include a type of mathematical optimization that has no closed-form solution and that is solved numerically, either online (trajectory optimization \cite{kuindersma_optimization-based_2016}, model predictive control \cite{schwenzer_review_2021}, etc.) or offline (reinforcement learning \cite{rudin_learning_2022}). Numerical tools, however, have a major drawback compared to simpler analytical approaches: the parameters of the problem do not appear explicitly in the solutions, which makes it much harder to generalize and reuse the results. Analytical solutions to control problems have the useful property of allowing the solution to be adjusted to different system parameters by simply substituting the new values in the equation. For instance, an analytical feedback law solution to a robot motion control problem can be transferred to a similar system by adjusting the values of parameters (lengths, masses, etc.) in the equation. However, with a reinforcement learning solution, we would have to re-conduct all the training, implying (generally) multiple hours of data collection and/or computation. It would be a great asset to have the ability to adjust black box numerical solutions with respect to some problem parameters.

%%%%%%%%%%%%%%%%%%%%%%
\begin{figure}[t]
    \begin{center}
        \includegraphics[width=8cm]{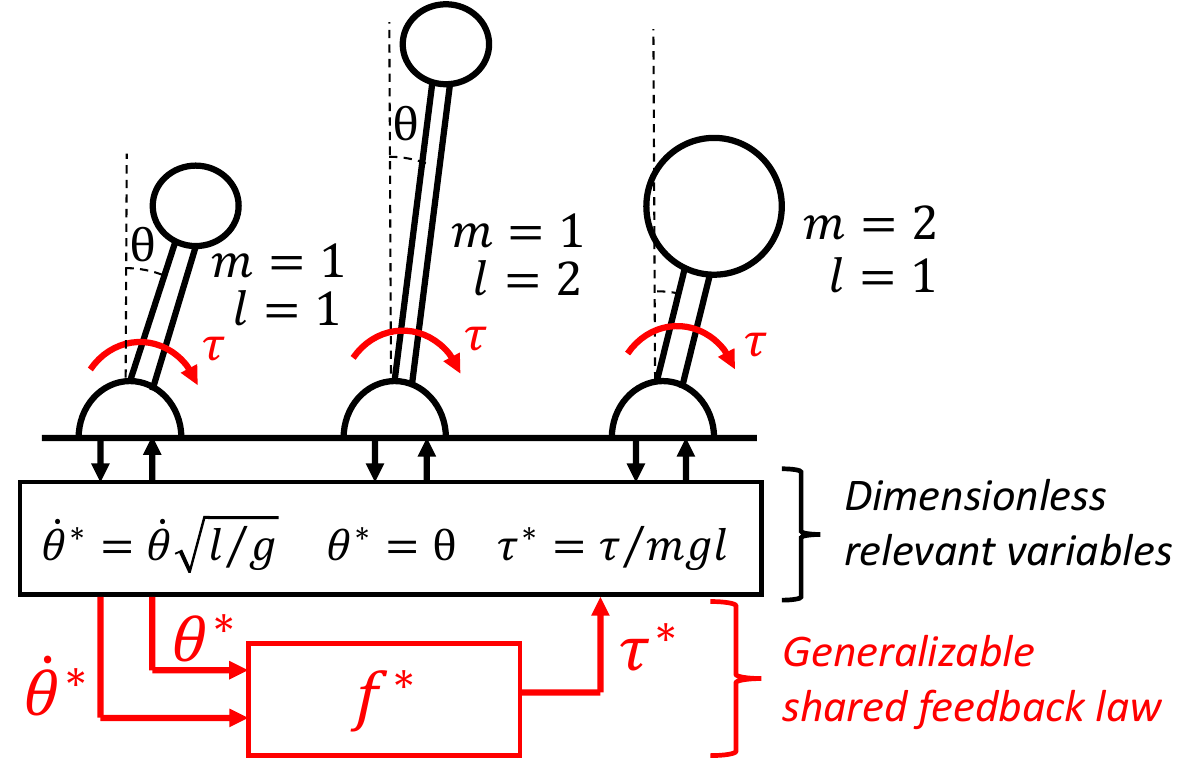}
        \caption{Shared dimensionless policy for inverted pendulums: Under some conditions various dynamic systems will share the same optimal policy up to scaling factors that can be found based on a dimensional analysis.}
        \label{fig:big_picture}
        \vspace{-10pt}
    \end{center}
\end{figure}
%%%%%%%%%%%%%%%%%%%%%%

In this paper, we explore the concept of dimensionless policies, a more generic form of knowledge conceptually illustrated at Figure \ref{fig:big_picture}, as an approach to generalize numerical solutions to motion control problems. First, in Section \ref{sec:dimenanalysis}, we use dimensional analysis (i.e., the Buckingham $\pi$ theorem \cite{buckingham_physically_1914}) to show that motion control problems with dimensionally similar context variables must share the same feedback law solution when expressed in a dimensionless form, and discuss the implications. Two main generic theoretical results, relevant for any physically meaningful control policies, are presented as Theorem \ref{theo:pistar} and Theorem \ref{theo:ab}. Then in Section \ref{sec:numcasestud} we present two case studies with numerical results. Optimal feedback laws computed with a dynamic programming algorithm are used to demonstrate the theoretical results and their relevance for 1) the classical motion control problem of swinging up an inverted pendulum in Section \ref{sec:optimalswingup} and 2) a car motion control problem in Section \ref{sec:optimalcar}. Furthermore, in Section \ref{sec:closedfrom}, we illustrate—with two examples—how the proposed dimensional scaling is equivalent to changing parameters in an analytical solution. 

A very promising application of the concept of dimensionless policies is to empower reinforcement learning schemes, for which data efficiency is critical \cite{sutton_reinforcement_2018}. For instance, it would be interesting to use the data of all vehicles on the road, even if they are of varying dimensions and dynamic characteristics, to learn appropriate maneuvers in situations that occur very rarely. This idea of reusing data or results in a different context is usually called transfer learning \cite{taylor_transfer_2009} and has received a great deal of research attention, mostly targeted at applying a learned policy to new tasks. The more specific idea of transferring policies and data between systems/robots has also been explored, with schemes based on modular blocks \cite{devin_learning_2017}, invariant features \cite{gupta_learning_2017}, a dynamic map \cite{helwa_multi-robot_2017}, a vector representation of each robot hardware \cite{chen_hardware_2018}, or using tools from adaptive control \cite{pereida_data-efficient_2018} and robust control \cite{sorocky_experience_2020}. Dimensionless numbers and dimensional analysis comprise a technique based on the idea that some relationships should not depend on units that can be used for analyzing many physical problems \cite{bertrand_sur_1878} \cite{rayleigh_viii_1892} \cite{buckingham_physically_1914}. The most well-known application in the field of fluid mechanics is the idea of matching ratios (i.e., Reynolds, Prandtl, or Mach numbers) to allow for the generalization of experimental results between systems of various scales. The recent success of machine learning and data-driven schemes bring front and center the question of generalizing results, and there is a renewed interest in using dimensional analysis in the context of learning \cite{bakarji_dimensionally_2022} \cite{fukami_robust_2021} \cite{xie_data-driven_2022}. In this paper, we present an initial exploration of how dimensional analysis can be applied specifically to help generalize policy solutions for motion control problems involving physically meaningful variables like force, length, mass, and time.

%%%%%%%%%%%%%%%%%%%%%%%%%%%%%%%%%%%%%%%%%%%%
\section{Dimensionless Policies}
\label{sec:dimenanalysis}

In the following section, we develop the concept of dimensionless policies based on the Buckingham $\pi$ theorem and present generic theoretical results that are relevant for any type of physically meaningful control policies. 

%%%%%%%%%%%%%%%%%%%%%%%%%%%%%%%%%%%%%%%%%%%%
\subsection{Context variables in the policy mapping}

\begin{figure}[htb]
    \vspace{-5pt}
    \begin{center}
        \vspace{-10pt}
        \subfloat[Generic policy]{\includegraphics[width=0.47\linewidth]{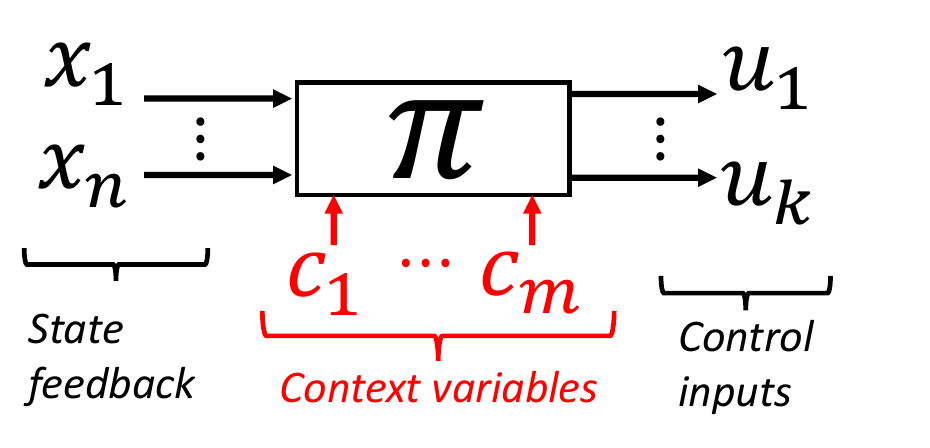}}
        \subfloat[Inverted pendulum example.]{\includegraphics[width=0.40\linewidth]{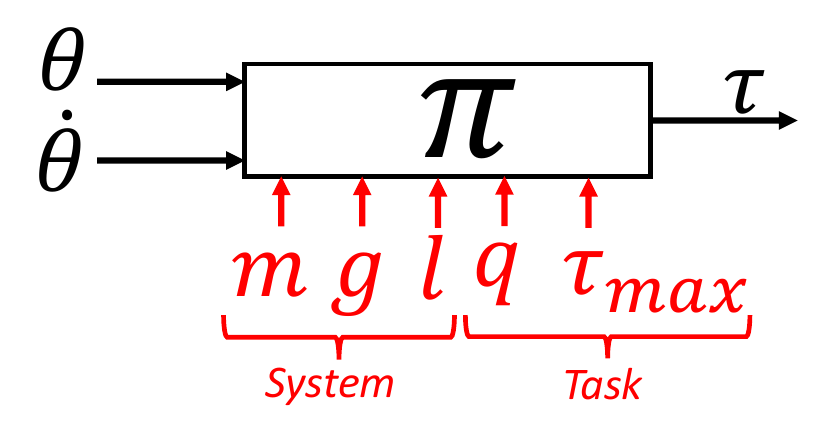}} 
        \caption{The policy $\pi$ is a feedback law that also includes problem parameters as additional arguments. 
        }\label{fig:policy_context}
    \end{center}
    \vspace{-5pt}
\end{figure}

Here, we call a feedback law a mapping $f$, specific to a given system, from a vector space representing the state $x$ of the dynamic system to a vector space representing the control inputs $u$ of the system:
%%%%%%%%%%%%%%%%%%%%%%
\begin{equation}
u
=
f \left(
x
\right)
\end{equation}
%%%%%%%%%%%%%%%%%%%%%%
Under some assumptions (fully observable systems, additive cost and infinite time horizon) the optimal feedback law is guarantee to be in this state feedback form \cite{bertsekas_dynamic_2012}. We will only consider motion control problems that lead to this type of time-independent feedback laws in the following analysis. To consider the question of how can this system-specific feedback law be transferred to a different context, it is useful to think about a higher dimension mapping $\pi$, which is herein referred to as a policy, also having  a vector of variables $c$ describing the context as an additional input argument as  illustrated in Figure \ref{fig:policy_context}. 

\begin{Definition}  
\label{def:policy}
A policy is defined as the solution to a motion control problem in the form of a function computing the control inputs $u$ from the system states $x$ and context parameters $c$ as follow:
%%%%%%%%%%%%%%%%%%%%%%
\begin{align}
u
=
&\pi \left(
x,
c
\right)
\label{eq:policy}\\
\textit{with} 
\quad 
u \in \mathbf{R}^k 
\quad 
&x \in \mathbf{R}^n 
\quad 
c \in \mathbf{R}^m 
\end{align}
%%%%%%%%%%%%%%%%%%%%%%
where $k$ is the dimension of the control input vector, $n$ is the dimension of the dynamic system state vector and $k$ is the dimension of the vector of context parameters.
\end{Definition}

The context $c$ is a vector of relevant parameters defining the motion control problem, i.e., parameters that affect the feedback law solution. The policy $\pi$ is thus a mapping consisting of the feedback law solutions for all possible contexts. In Section \ref{sec:optimalswingup}, a case study is conducted by considering the optimal feedback law for swinging up a torque-limited inverted pendulum. For this example, the context variables are the pendulum mass $m$, the gravitational constant $g$, and the length $l$, as well as what we call task parameters: a weight parameter in the cost function $q$ and a constraint $\tau_{max}$ on the maximum input torque. For a given pendulum state, the optimal torque is also a function of the context variables, i.e., the solution is different if the pendulum is heavier or more torque limited.
%%%%%%%%%%%%%%%%%%%%%%
%%%%%%%%%%%%%%%%%%%%%%

\begin{figure}[htb]
    \begin{center}
        \includegraphics[width=0.95\linewidth]{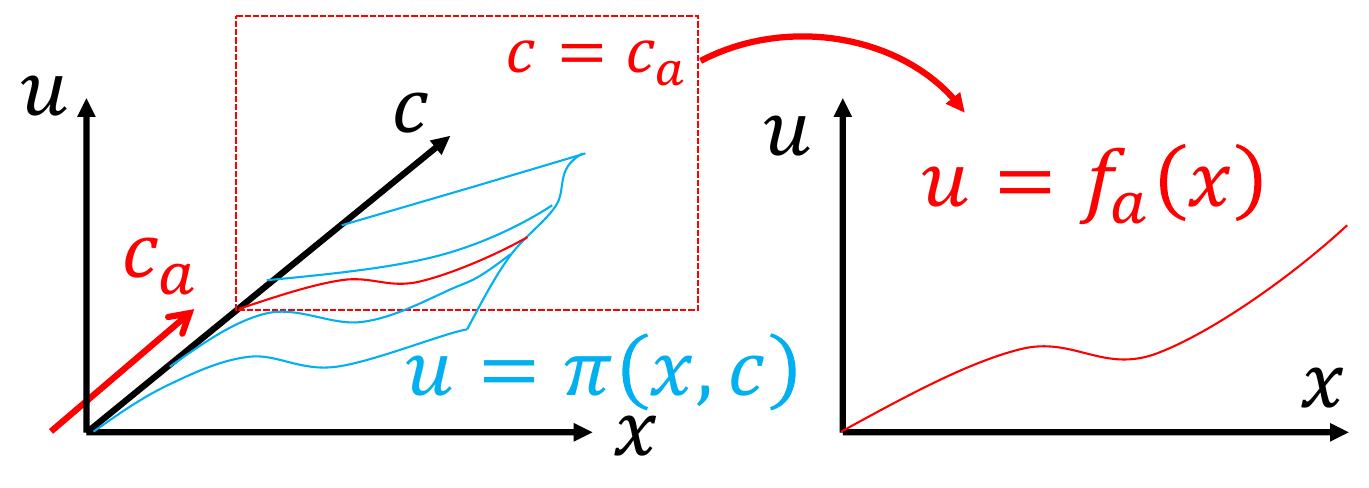}
        \caption{A feedback law $f$ is a slice of the higher dimensional policy mapping $\pi$ in a specific context.
        }\label{fig:slice}
    \end{center}
\end{figure}

\begin{Definition}  
\label{def:fa}
A feedback law with a subscript letter $a$ is defined as the solution to a motion problem for a specific situation defined by an instance of context variables $c_a$, as follow:
%%%%%%%%%%%%%%%%%%%%%%
\begin{equation}
f_a \left(
x 
\right) 
=
\pi \left(
x,
c = c_a
\right) \quad \quad \forall x
\end{equation}
%%%%%%%%%%%%%%%%%%%%%%
\end{Definition}
The feedback law $f_a$ thus represents a slice of the global policy when the context variables are fixed at $c_a$ values as illustrated in Figure \ref{fig:slice}.  
%%%%%%%%%%%%%%%%%%%%%%
%%%%%%%%%%%%%%%%%%%%%%

The goal of generalizing a feedback law to a different context can thus be formulated into the following question: If a feedback law $f_a$ is known for a context described by variables $c_a$, can this knowledge help us deduce the policy solution in a different context, namely $c_b$?
%%%%%%%%%%%%%%%%%%%%%%
\begin{equation}
\pi \left(
x,
c = c_a
\right) = 
f_a \left(
x 
\right) 
\quad \Rightarrow \quad
\pi \left(
x,
c = c_b
\right) = \, ? \quad
\end{equation}
%%%%%%%%%%%%%%%%%%%%%%
Using the Buckingham $\pi$ theorem \cite{buckingham_physically_1914}, we will show that if the context is dimensionally similar, then both feedback laws must be equal up to scaling factors (Theorem \ref{theo:ab}).
% %%%%%%%%%%%%%%%%%%%%%%
% \begin{equation}
% \underbrace{\begin{bmatrix}
% u_1 \\
% \vdots \\
% u_k
% \end{bmatrix}}_{\text{inputs}}
% =
% \pi \left(
% \underbrace{\begin{bmatrix}
% x_1 \\
% \vdots \\
% x_n
% \end{bmatrix}}_{\text{states}}
% ,
% \underbrace{
% \underbrace{\begin{bmatrix}
% c_1 \\
% \vdots \\
% \vdots 
% \end{bmatrix}}_{\text{system}}
% ,
% \underbrace{\begin{bmatrix}
% \vdots \\
% c_{m}
% \end{bmatrix}}_{\text{task}}
% }_{\text{Context $c$}}
% \right) 
% \label{eq:vectorpolicy}
% \end{equation}
% %%%%%%%%%%%%%%%%%%%%%%

\subsection{ Buckingham $\pi$ theorem  }
The Buckingham $\pi$ theorem \cite{buckingham_physically_1914}, is a tool based on dimensional analysis \cite{bertrand_sur_1878} \cite{rayleigh_viii_1892}, that allow to restate a relationship involving multiple physically meaningful dimensional variables, using a lesser number of dimensionless variables:
\begin{equation}
x_1 = f( x_2, \; \hdots \;, x_n )
\quad  \Rightarrow  \quad
\Pi_1 = f( \Pi_2, \hdots, \Pi_{p} )
\label{eq:bpi}
\end{equation}
If $d$ fundamental dimensions are involved in the $n$ dimensional variables (for instance time [T], length [L] and mass [M]), then the number of required dimensionless variables, often called $\Pi$ groups is
$p \geq n - d$. In most situations, the number of variables in the relationship can be reduced directly by the number of fundamental dimensions involved and $p = n - d$. The Buckingham $\pi$ theorem provides a methodology to generate the $\Pi$ groups, however the choice of $\Pi$ groups is not unique. The approach is to select (arbitrarily) $d$ variables involving the $d$ fundamental dimensions independently, called the repeated variables. Then, the $\Pi$ groups are generated by multiplicating all the other variables, by the repeated variables exponentiated by rational exponents selected to make the group dimensionless. Assuming $x_1$, ... $x_d$, where the selected repeated variables, the $\Pi$ groups are:
\begin{equation}
\Pi_i = x_{d+i} \; 
\underbrace{
x_1^{e_{1i}} \; x_2^{e_{2i}} \hdots x_d^{e_{di}} 
}_{\text{Repeated variables}}
\quad i=\{1, ... , p \}
\label{eq:bpi}
\end{equation}
Finding the correct exponents to make all group dimensionless can be formulated as solving a linear system of $d$ equations. We refer to previous literature for more details on the theorem, and here use it specifically on the defined concept of policy map.

\subsection{Dimensional analysis of the policy mapping}
\label{sec:dimpolicymap}

If a policy is physically meaningful (for example, a policy that computes a force based on position and velocity, but not a policy for playing chess), we can use the Buckingham $\pi$ theorem to simplify the policy in dimensionless form.

\begin{Theorem}\label{theo:pistar}
If a policy is physically meaningful and all its variables involve $d$ fundamental dimensions that are independently present in the context variables $c$, then the policy can be restated in a dimensionless form as follow:
%%%%%%%%%%%%%%%%%%%%%%
\begin{align}
u = \pi( x , c )
\quad \quad  &\Rightarrow \quad 
u^* = \pi^*( x^* , c^* )
\label{eq:vectordimpolicyshort}
\\
u \in \mathbf{R}^k 
\;
x \in \mathbf{R}^n 
\; 
c \in \mathbf{R}^{m}
\;
\;
&
\;
\;
u^* \in \mathbf{R}^k 
\;
x^* \in \mathbf{R}^n 
\;
c^* \in \mathbf{R}^{(m-d)}
\end{align}
%%%%%%%%%%%%%%%%%%%%%%
where the dimensionless variables can be related to dimensional variables using transformation matrices that depends only on the context variables as follow:
%%%%%%%%%%%%%%%%%%%%%%
\begin{align}
u^* &= \left[ T_u(c) \right] \, u  \label{eq:Tu} \\
x^* &= \left[ T_x(c) \right] \, x \label{eq:Tx} \\
c^* &= \left[ T_c(c) \right] \, c \label{eq:Tc}
\end{align}
%%%%%%%%%%%%%%%%%%%%%%
Furthermore, the transformation matrices can be used to relate the dimensional and dimensionless policy as follow:
%%%%%%%%%%%%%%%%%%%%%%
\begin{align}
\pi( x , c ) = 
T_u^{-1}(c) \; 
\pi^* \Big( \; 
T_x(c) x 
\; , \; 
T_c(c) c  
\; \Big)
\label{eq:pistar2pi}
\end{align}
%%%%%%%%%%%%%%%%%%%%%%
\end{Theorem}

\begin{proof}
For a system with $k$ control inputs, we can treat the policy as $k$ mappings from states and context variables to each scalar control input $u_j$:
%%%%%%%%%%%%%%%%%%%%%%
\begin{equation}
u_j = \pi_j \left(
x_1, \hdots, x_n, 
c_1, \hdots \hdots, c_{m}
\right) 
\label{eq:scalarpolicy}
\end{equation}
%%%%%%%%%%%%%%%%%%%%%%
where Equation \eqref{eq:scalarpolicy} is the $j$th line of the policy in vector form, as described by Equation \eqref{eq:policy}. Then, if the state vector is defined by $n$ variables, and the context is defined by $m$ (system and task) parameters, then each mapping $\pi_j$ is a relation between $1 + n + m $ variables. Under the assumption that the policy involves physically meaningful variables, and that it is invariant under an arbitrary scaling of any fundamental dimensions-- i.e. independent of a system of units--, then we can apply the Buckingham $\pi$ theorem \cite{buckingham_physically_1914} to conclude that if $d$ dimensions are involved in all of those variables, then Equation \eqref{eq:scalarpolicy} can be restated into an equivalent relationship between $p$ dimensionless $\Pi$ groups where $p \geq 1 + n + m  - d$. Assuming that $d$ dimensions are involved in the $m$ context variables, and that we are in the typical scenario where maximum reduction is possible  ($p = 1 + n + m  - d$), then we can select $d$ context variables $\{c_1, c_2 , \hdots,c_d\}$ as the basis (the repeated variables) to scale all other variables into dimensionless $\Pi$ groups. We denote dimensionless $\Pi$ group as the base variables with an asterisk (*), as follows:
%%%%%%%%%%%%%%%%%%%%%%
\begin{align}
u_j^* &= u_j \left[ c_1 \right]^{e^u_{1j}} \left[ c_2 \right]^{e^u_{2j}} \hdots \left[ c_d \right]^{e^u_{dj}} \quad \scriptstyle j = \{ 1, \hdots , k \} \label{eq:piu}\\
x_i^* &= x_i \left[ c_1 \right]^{e^x_{1i}} \left[ c_2 \right]^{e^x_{2i}} \hdots \left[ c_d \right]^{e^x_{di}} \quad \scriptstyle i = \{ 1, \hdots , n \} \label{eq:pix}\\
c_i^* &= c_i\left[ c_1 \right]^{e^c_{1i}} \left[ c_2 \right]^{e^c_{2i}} \hdots \left[ c_d \right]^{e^c_{di}} \quad \scriptstyle i = \{ d+1, \hdots , m \} \label{eq:pic}%\\
%b_i^* &= b_i \left[ a_1 \right]^{e_1} \left[ a_2 \right]^{e_2} \hdots \left[ a_d \right]^{e_d} \quad  i = \{ 1, \hdots , l \} 
\end{align}
%%%%%%%%%%%%%%%%%%%%%%
where exponents $e_{ij}$ are rational numbers selected to make all equations dimensionless. We can then define transformation matrices and write Equations \eqref{eq:piu}, \eqref{eq:pix}, and \eqref{eq:pic} in a vector form where the repeated variables are grouped into matrices defined as shown at Equations \eqref{eq:TU},  \eqref{eq:TX} and  \eqref{eq:TC}
\begin{figure*}[!t]
% ensure that we have normalsize text
\normalsize
\setcounter{equation}{17}
%%%%%%%%%%%%%%%%%%%%%%
\begin{align}
\underbrace{
\scriptsize
\begin{bmatrix}
u_1^* \\ \vdots \\ u_k^*
\end{bmatrix}
}_{u^*}
&= 
\underbrace{
\scriptsize
\begin{bmatrix}
\left( \left[ c_1 \right]^{e^u_{11}} \left[ c_2 \right]^{e^u_{21}} \hdots \left[ c_d \right]^{e^u_{d1}} \right ) & 0 & 0  \\
0 & \ddots & 0 \\
0 & 0 & \left( \left[ c_1 \right]^{e^u_{1k}} \left[ c_2 \right]^{e^u_{2k}} \hdots \left[ c_d \right]^{e^u_{dk}} \right)
\end{bmatrix}
}_{T_u(c)}
\underbrace{
\scriptsize
\begin{bmatrix}
u_1 \\ \vdots \\ u_k
\end{bmatrix}
}_{u}
%%%%%%%%%%%%%%%%%%%%%%
\label{eq:TU}
\\ 
%%%%%%%%%%%%%%%%%%%%%%
\underbrace{
\scriptsize
\begin{bmatrix}
x_1^* \\ \vdots \\ x_n^*
\end{bmatrix}
}_{x^*}
&=
\underbrace{
\scriptsize
\begin{bmatrix}
\left( \left[ c_1 \right]^{e^x_{11}} \left[ c_2 \right]^{e^x_{21}} \hdots \left[ c_d \right]^{e^x_{d1}} \right ) & 0 & 0  \\
0 & \ddots & 0 \\
0 & 0 & \left( \left[ c_1 \right]^{e^x_{1n}} \left[ c_2 \right]^{e^x_{2n}} \hdots \left[ c_d \right]^{e^x_{dn}} \right)
\end{bmatrix}
}_{T_x(c)}
\underbrace{
\scriptsize
\begin{bmatrix}
x_1 \\ \vdots \\ x_n
\end{bmatrix}
}_{x}
%%%%%%%%%%%%%%%%%%%%%%
\label{eq:TX}
\\ 
%%%%%%%%%%%%%%%%%%%%%%
\underbrace{
\scriptsize
\begin{bmatrix}
c_{d+1}^* \\ \vdots \\ c_m^*
\end{bmatrix}
}_{c^*}
&=
\underbrace{
\scriptsize
\begin{bmatrix}
0 & \hdots & 0 &  \left( \left[ c_1 \right]^{e^u_{1(d+1)}} \left[ c_2 \right]^{e^u_{2(d+1)}} \hdots \left[ c_d \right]^{e^u_{d(d+1)}} \right ) & 0 & 0  \\
0 & \hdots & 0 & 0 & \ddots & 0 \\
%\underbrace{
0 & \hdots & 0 
%}_{\text{dim=d}}
& 0 & 0 & \left( \left[ c_1 \right]^{e^u_{1m}} \left[ c_2 \right]^{e^u_{2m}} \hdots \left[ c_d \right]^{e^u_{dm}} \right)
\end{bmatrix}
}_{T_c(c)}
\underbrace{
\scriptsize
\begin{bmatrix}
c_1 \\ \vdots \\ c_m
\end{bmatrix}
}_{c}
\label{eq:TC}
\end{align}
%%%%%%%%%%%%%%%%%%%%%%
%\setcounter{equation}{\value{MYtempeqncnt}}
% The spacer can be tweaked to stop underfull vboxes.
\vspace*{4pt}
\end{figure*}
%%%%%%%%%%%%%%%%%%%%%%
which correspond to Equations \eqref{eq:Tu}, \eqref{eq:Tx} and \eqref{eq:Tc}. Matrices $T_u$ and $T_x$ are square diagonal matrices and Equations \eqref{eq:Tu} and \eqref{eq:Tx} are thus inversibles (unless a repeated variable is equal to zero) and can be used to go back and forth between dimensional and dimensionless states and input variables. Matrix $T_c$ consist in a block of $d$ columns of zeros, followed by a diagonal block of dimensions $(m-d) \times (m-d)$, and Equation \eqref{eq:Tc} is not inversible. For a given context $c$, there is only one dimensionless context $c^*$, however a given dimensionless context $c^*$ may correspond to multiple dimensional contexts $c$. 

Then, the Buckingham $\pi$ theorem tell us that the relationship described by Equation \eqref{eq:scalarpolicy} can be restated in a relationship between the $\Pi$ groups involving $d$ less variables, which based on the selected repeated variable correspond to:
%%%%%%%%%%%%%%%%%%%%%%
\begin{equation}
u_j^* = \pi_j^* \left(
x_1^*, \hdots, x_n^*, 
c_{d+1}^*, \hdots, c_{m}^*
\right) 
\label{eq:scalardimpolicy}
\end{equation}
%%%%%%%%%%%%%%%%%%%%%%
By applying the same procedure to all control inputs, we can then assemble all $k$ mappings back into a vector form, as follows:
%%%%%%%%%%%%%%%%%%%%%%
\begin{equation}
\underbrace{
\begin{bmatrix}
u_1^* \\
\vdots \\
u_k^*
\end{bmatrix}
=
\pi^* \Biggl(
\begin{bmatrix}
x_1^* \\
\vdots \\
x_n^*
\end{bmatrix}
}_{\text{Dimensionless feedback law $f^*$}}
,
\underbrace{
\begin{bmatrix}
c_{d+1}^* \\
\vdots \\
c_{m}^*
\end{bmatrix}
% ,
% \begin{bmatrix}
% c_{m+1}^* \\
% \vdots \\
% c_{m+l}^*
% \end{bmatrix}
}_{\text{context $c^*$}}
\Biggr)
\label{eq:vectordimpolicy}
\end{equation}
%%%%%%%%%%%%%%%%%%%%%%
that correspond to Equation \eqref{eq:vectordimpolicyshort}. Finally, based on the defined transformations at Equations \eqref{eq:Tu}, \eqref{eq:Tx} and \eqref{eq:Tc} we can relate the dimensional policy to the dimensionless version as follow:
%%%%%%%%%%%%%%%%%%%%%%
\begin{align}
\pi( x , c ) = 
\underbrace{
T_u^{-1}(c) \; 
\underbrace{
\pi^* \Big( \; 
\underbrace{
T_x(c) x 
}_{x^*}
\; , \; 
\underbrace{
T_c(c) c  
}_{c^*}
\; \Big)
}_{u^*}
}_{u}
\end{align}
%%%%%%%%%%%%%%%%%%%%%%
which correspond to Equation \eqref{eq:pistar2pi}.
\end{proof}

\subsection{Transferring feedback laws between similar systems}\label{sec:transfer}
Based on the dimensional analysis, we can demonstrate that any feedback law can be generalized to a different context, under the condition of dimensional similarity. In this section, we show that a feedback law can be transferred exactly to another motion control problem by scaling the input and output of the function based on matrices that can be computed using the dimensional analysis. The salient feature of this result is that the conditions are very generic, even a black-box discontinuous non-linear policy (such as those obtained using deep-reinforcement learning algorithms) can be transferred this way. The limitation is that the condition for an exact transfer is having equal dimensionless context variables $c^*$.

First, it is useful to define dimensionless feedback laws that correspond to specific cases of the dimensionless policy, as we defined for the dimensional mapping.
\begin{Definition}  \label{def:fastar}
We denote a dimensionless feedback law $f_a^*$, the global dimensionless policy for a specific instance of context variables $c_a$, as follow:
% \begin{Definition}
% A dimensionless feedback law noted $f_a^*$ correspond to
%%%%%%%%%%%%%%%%%%%%%%
\begin{equation}
f_a^*( x^* ) = \pi^*( x^* , c^* = c^*_a ) \quad \forall x^*
\label{eq:dimfeedback}
\end{equation}
% %%%%%%%%%%%%%%%%%%%%%%
where $c^*_a$ is the dimensionless version of the context variables instance $c_a$, and equal to:
%%%%%%%%%%%%%%%%%%%%%%
\begin{equation}
c^*_a = T_c(c_a) \, c_a
\end{equation}
% %%%%%%%%%%%%%%%%%%%%%%
\end{Definition}

\begin{Lemma}\label{lemma:star2star}
Two feedback laws, that are solutions to the same motion control problem for two instance of context variables, will be equal in dimensionless form if they share the same dimensionless context:
%%%%%%%%%%%%%%%%%%%%%%
\begin{equation}
f_b^*(x^*) = f_a^*(x^*) \quad \forall x^*  \quad  \text{if} \quad c_a^* = c_b^* 
\label{eq:starstar}
\end{equation}
%%%%%%%%%%%%%%%%%%%%%%
\end{Lemma}
\begin{proof}
This follow from the definition:
%%%%%%%%%%%%%%%%%%%%%%
\begin{align}
f_a^*(x^*) &= \pi^*( x^* , c^* = c_a^*)  \\
f_a^*(x^*) &= \pi^*( x^* , c^* = c_b^*)  \\
f_a^*(x^*) &= f_b^*(x^*) 
\end{align}
%%%%%%%%%%%%%%%%%%%%%%
%The global policy $\pi^*$ is a function and give the same output for the same input.
\end{proof}
\begin{Lemma}\label{lemma:f2star}
In a specific context described by variables $c_a$, a dimensional feedback law can be restated into a dimensionless form, and vice versa, by scaling the input and the output using the defined transformation matrices $T_x(c_a)$ and $T_u(c_a)$ as follow:
%%%%%%%%%%%%%%%%%%%%%%
\begin{align}
f_a ( x ) &= T^{-1}_u(c_a) 
\underbrace{
f_a^* \Big(  
\underbrace{
T_x(c_a) \; x
}_{x^*}
\Big)
}_{u^*} \quad \forall x
\label{eq:star2f}
\\
f_a^* ( x^* ) &= T_u(c_a) 
\underbrace{
f_a \Big(  
\underbrace{
T_x^{-1}(c_a) \; x^*
}_{x}
\Big)
}_{u} \quad \forall x^*
\label{eq:f2star}
\end{align}
%%%%%%%%%%%%%%%%%%%%%%
\end{Lemma}
\begin{proof} 
Starting from Equation \ref{eq:pistar2pi} and substituting $c^*$ with a specific instance $c^*_a$, then substituting policy maps on each side with feedback laws $f_a$ and $f_a^*$ based on the definition, we obtain Equation \ref{eq:star2f}:
%%%%%%%%%%%%%%%%%%%%%%
\begin{align}
% \pi( x , c ) &= T_u^{-1}(c) \; 
% \pi^* \Big( \; 
% T_x(c) x 
% \; , \; 
% T_c(c) c  
% \; \Big)
% \quad \forall (x,c)
% \\
\pi( x , c_a ) &= T_u^{-1}(c_a) \; 
\pi^* \Big( \; 
T_x(c_a) x 
\; , \; 
T_c(c_a) c_a
\; \Big) 
\\
f_a ( x ) &= T_u^{-1}(c_a) f_a^* \Big( T_x(c_a) \, x \Big)
\end{align}
%%%%%%%%%%%%%%%%%%%%%%
Then, starting from the right side of Equation \ref{eq:f2star} and substituting the function $f_a$ with Equation $\ref{eq:star2f}$, the matrices are reduced to identity matrices and we obtain Equation \ref{eq:f2star}:
%%%%%%%%%%%%%%%%%%%%%%
\begin{align}
T_u(c_a) T_u^{-1}(c_a) 
f_a^* \Big(  
T_x^{-1}(c_a) T_x(c_a)  \; x^*
\Big)
&= f_a^*(x^*)
\end{align}
%%%%%%%%%%%%%%%%%%%%%%
% %%%%%%%%%%%%%%%%%%%%%%
% \begin{align}
% f_a ( x ) &= T^{-1}_u(c_a)  f_a^* \Big( T_x(c_a) \; x \Big) \quad \forall x \\
% f_a ( x ) &= T^{-1}_u(c_a)  T_u(c_a)  f_a \Big( T_x^{-1}(c_a) T_x(c_a) \; x \Big) \quad \forall x\\
% f_a ( x ) &= f_a ( x ) \quad \forall x
% \end{align}
% %%%%%%%%%%%%%%%%%%%%%%
\end{proof}
\begin{Theorem}\label{theo:ab}
If a feedback law $f_a$ is known—for instance, as the result of a numerical algorithm—and this is the solution to a motion control problem with context variables $c_a$, we can compute the solution $f_b$ to the same motion control problem for different context variables $c_b$ by scaling the input and output of $f_a$ as follow:
%%%%%%%%%%%%%%%%%%%%%%
\begin{align}
f_b ( x ) &= 
\left[ T^{-1}_u(c_b) 
T_u(c_a) \right] \,
f_a \left( 
\left[
T_x^{-1}(c_a) 
T_x(c_b)
\right] \,
x
\right) \quad \forall x
\label{eq:ab_transform}
\end{align}
%%%%%%%%%%%%%%%%%%%%%%
if the contexts $c_a$ and $c_b$ are dimensionally similar, i.e., if the following condition is true:
%%%%%%%%%%%%%%%%%%%%%%
\begin{equation}
T_c( c_b ) \; c_b  = T_c( c_a ) \; c_a
\label{eq:dimcontextequal}
\end{equation}
%%%%%%%%%%%%%%%%%%%%%%
\end{Theorem}
\begin{proof}
First, $f_b$ can be written based on its dimensionless form $f_b^*$ in a context $c_b$ using Equation \eqref{eq:star2f} from Lemma \ref{lemma:f2star}. Also, based on Lemma \ref{lemma:star2star}, under the similarity condition—i.e. $c_b^* = c_a^*$ or equivalently $T(c_b)c_b = T(c_a)c_a$— we have that $f^*_b$ is equal to $f^*_a$. Finally, $f^*_a$ can be written based on its dimensional form $f_a$ in a context $c_a$, using Equation \eqref{eq:f2star} from Lemma \ref{lemma:f2star}, as follow:
%%%%%%%%%%%%%%%%%%%%%%
\begin{align}
f_b ( x ) &=  T^{-1}_u(c_b)  f_b^* \Big( T_x(c_b) \; x \Big) 
\\
f_b ( x ) &=  T^{-1}_u(c_b)  f_a^* \Big( T_x(c_b) \; x \Big) \quad \quad \text{if} \quad c_b^* = c_a^*
\\
f_b ( x ) &= 
\left[ T^{-1}_u(c_b) 
T_u(c_a) \right] \,
f_a \left( 
\left[
T_x^{-1}(c_a) 
T_x(c_b)
\right] \,
x
\right)
\end{align}
%%%%%%%%%%%%%%%%%%%%%%
\end{proof}

\begin{figure}[htb]
\begin{center}
\includegraphics[width=0.45\textwidth]{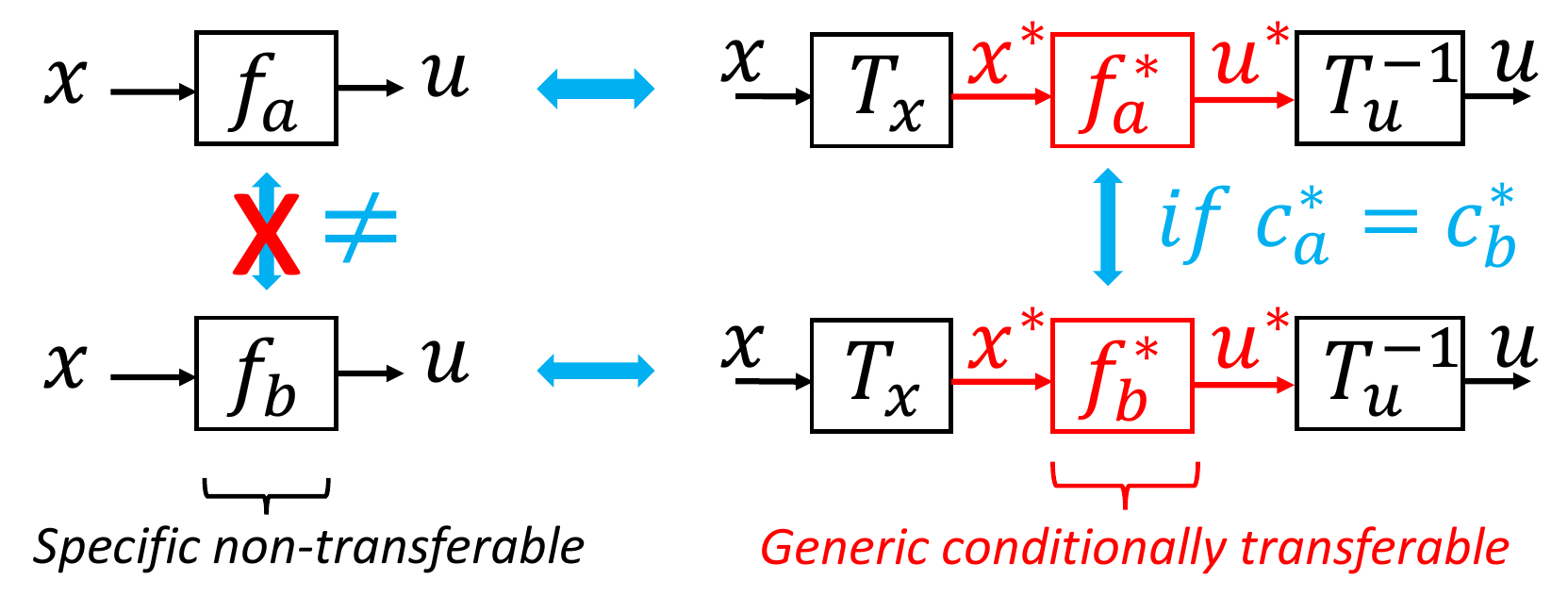}
\caption{Isolating the dimensionless knowledge in a policy allow its exact transfer to any dimensionally similar motion control problem.}\label{fig:dimpol}
\end{center}
\end{figure}

The idea is summarized in Figure \ref{fig:dimpol}. To transfer a feedback law, we must first extract the dimensionless form, a more generic form of knowledge, and then scale it back to the new context. 
%%%%%%%%%%%%%%%%%%%%%%
%%%%%%%%%%%%%%%%%%%%%%

\subsection{Dimensionally similar contexts}

\begin{figure}[ht]
\begin{center}
\includegraphics[width=0.95\linewidth]{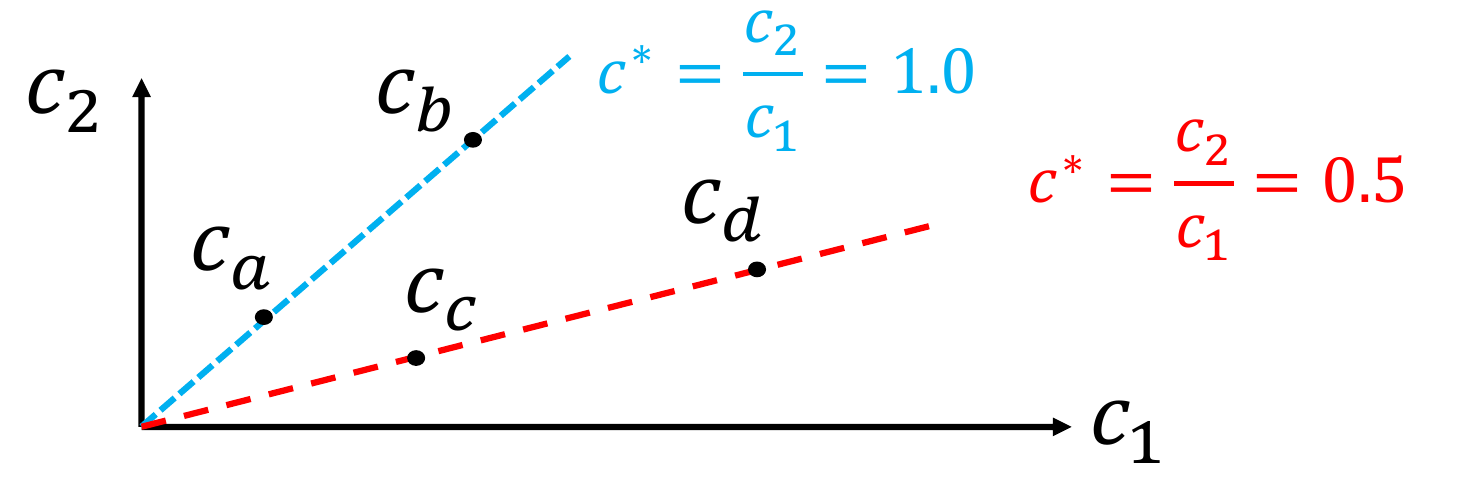}
\caption{Example of dimensionally similar contexts subsets that are lines in a plane ($m=2$ and $d=1$). Context $c_a$ is dimensionally similar to $c_b$ but not to $c_c$ or $c_d$.}\label{fig:c_space}
\end{center}
\end{figure}

\begin{figure}[ht]
\begin{center}
\includegraphics[width=0.95\linewidth]{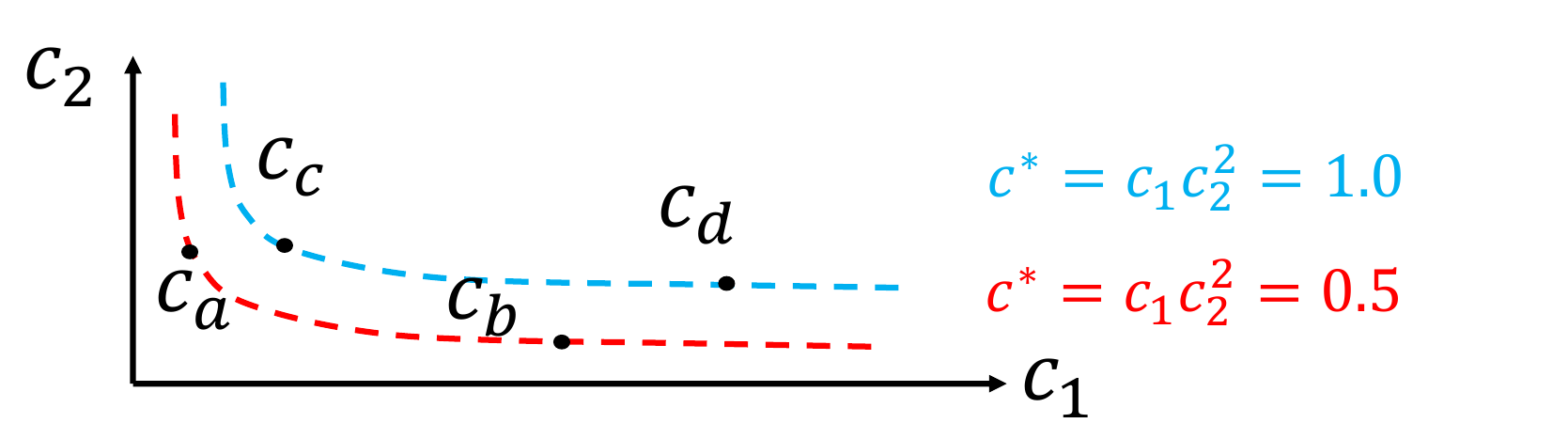}
\caption{Example of dimensionally similar contexts subsets that are non-linear curves in a plane ($m=2$ and $d=1$). Context $c_a$ is dimensionally similar to $c_b$ but not to $c_c$ or $c_d$.}\label{fig:c_space2}
\end{center}
\end{figure}

Equation \eqref{eq:ab_transform} can be used to scale a policy for an exact transfer of policy solutions between context $c$ sharing the same dimensionless context $c^*$, a condition that is refer to as \textit{dimensionally similar}. Equation \eqref{eq:Tc} is a mapping from a $m$ dimensional space to a $m-d$ dimensional space, and its inverse has multiple solutions. A given dimensionless context $c^*$ corresponds to a subset of all possible values of dimensional context $c$. As illustrated at Figure \ref{fig:c_space} and  Figure \ref{fig:c_space2} with low-dimensional examples ($m=2$ and $d=1$), the subsets of context $c$ leading to the same $c^*$ can be linear if $c^*$ is just a ratio of two variables of the same dimension, or a non-linear curve if $c^*$ involves exponents leading to more a complex polynomial relationship. In general when the context $c$ involves many dimensions, it is important to note that the similarity condition means meeting multiple conditions (one for each element of the vector $c^*$) in a higher-dimensional space as illustrated at Figure \ref{fig:context_3d} for the pendulum swing-up example that is studied in the next section. To some degree, this dimensionally similar context condition is a technique to regroup the motion control problems that are the same up to scaling factors. Therefore, it is also logical that their solutions should be equivalent up to scaling factors. 

%%%%%%%%%%%%%%%%%%%%%%
%%%%%%%%%%%%%%%%%%%%%%

%%%%%%%%%%%%%%%%%%%%%%
%%%%%%%%%%%%%%%%%%%%%%

\subsection{Summary of the theoretical results}

The dimensional analysis lead us to the following relevant theoretical results, that are very generic since no assumptions on the form of the policy function are necessary: 
\begin{enumerate}
    \item The global problem of learning $\pi(x,c)$, i.e., the feedback policies for all possible contexts, is simplified in a dimensionless form $\pi^*(x^*,c^*)$ because we can remove $d$ input dimensions from the unknown mapping (typically, $d$ would be 2 or 3 for controlling a physical system involving time, force, and length), see Theorem \ref{theo:pistar}. 
    \item The feedback law solutions of dimensionally similar subset of contexts share the exact same solution when restated in a dimensionless form, see Lemma \ref{lemma:star2star}. 
    \item A feedback law, which is a solution to a motion control problem in a context can be transferred exactly to another context, under a condition of dimensional similarity, by scaling appropriately its inputs and outputs, see Theorem \ref{theo:ab}.
\end{enumerate}

Just for illustrating purposes, lets imagine we have a policy for a spherical submarine where the context is defined by a velocity, a viscosity and a radius. In dimensionless form we would find that the context can be described by a single variable, the Reynolds number, and that \textbf{1)} learning the policy will be easier in dimensionless form because it is a function of a lesser number of variables and \textbf{2)} that if we know the feedback law solution for a specific context of velocity, viscosity and radius, then we can actually re-use it for multiple versions of the same motion control problem sharing the same Reynolds number.

\section{Case studies with numerical results}
\label{sec:numcasestud}
In this section, we use numerically generated optimal policy solutions for two motion control problem as examples illustrating the salient features of the presented theoretical results of section \ref{sec:dimenanalysis} and the potential for transfer learning.

%%%%%%%%%%%%%%%%%%%%%%
\subsection{Optimal pendulum swing-up task}
\label{sec:optimalswingup}
The first numerical example is the classical pendulum swing-up task. This example illustrates that an optimal feedback law in the form of a table look-up generated for a pendulum of a given mass and length, can be transferred to a pendulum of a different mass and length if the motion control problem is dimensionally similar. The example is also used to introduce the concept of regime for motion control problem.

\subsubsection{Motion control problem}
The motion control problem is defined here as finding a feedback law for controlling the dynamic system described by the following differential equation:
%%%%%%%%%%%%%%%%%%%%%%
\begin{equation}
ml^2 \ddot{\theta} - mgl \sin \theta = \tau
\label{eq:pendulum_dynamics}
\end{equation}
%%%%%%%%%%%%%%%%%%%%%%
which minimizes the infinite horizon quadratic cost function given by:
%%%%%%%%%%%%%%%%%%%%%
\begin{equation}
J %= \int_0^{\infty}{( q^2 \theta^2 + 0 \, \dot{\theta}^2 + 1 \, \tau^2 ) dt }
% = \int_0^{\infty}{\left(
% \begin{bmatrix}
% \theta & \dot{\theta}
% \end{bmatrix}^T
% \underbrace{
% \begin{bmatrix}
% q^2 &  0 \\ 0 & 0
% \end{bmatrix}}_{Q}
% \begin{bmatrix}
% \theta \\ \dot{\theta}
% \end{bmatrix}
%  + 
%  \begin{bmatrix}
% \tau
% \end{bmatrix}^T
% \underbrace{
% \begin{bmatrix}
% 1
% \end{bmatrix}}_{R}
% \begin{bmatrix}
% \tau
% \end{bmatrix}
% \right) dt } 
= \int_0^{\infty}{ \left( q^2 \theta^2 \, + \, \tau^2 \right) dt }
\label{eq:pendulum_cost}
\end{equation}
%%%%%%%%%%%%%%%%%%%%%%
subject to input constraints given by:
%%%%%%%%%%%%%%%%%%%%%%
\begin{equation}
- \tau_{max} \leq \tau \leq \tau_{max}
\label{eq:pendulum_constraints}
\end{equation}
%%%%%%%%%%%%%%%%%%%%%%
Note that, here, \textbf{1)} the cost function parameter $q$ has a power of two to allow its value to be in units of torque; \textbf{2)} it was chosen not to penalize high velocity values for simplicity; \textbf{3)} the weight multiplying the torque is set to one without a loss of generality, as only the relative values of weights impact the optimal solution; and \textbf{4)} all parameters are time-independent constants. Thus, assuming that there is no hidden variables and that Equations \eqref{eq:pendulum_dynamics}, \eqref{eq:pendulum_cost}, and \eqref{eq:pendulum_constraints} fully describe the problem, the solution—i.e., the optimal policy for all contexts—involves the variables listed in Table \ref{tb:optimalswingup}, and should be of the form given by:
%%%%%%%%%%%%%%%%%%%%%%
\begin{equation}
\underbrace{\tau}_{\text{inputs}}
=
\pi \left(
\underbrace{ \theta, \dot{\theta} }_{\text{states}},
\underbrace{
\underbrace{ m , g , l }_{\text{system parameters}},
\underbrace{ q , \tau_{max} }_{\text{task parameters}}
}_{\text{Context $c$}}
\right)
\label{eq:pendulumpi}
\end{equation}
%%%%%%%%%%%%%%%%%%%%%%

%%%%%%%%%%%%%%%%%%%%%%%%%%%%%%%%%%%%%%%%%%%%
\begin{table}[htb]
    \centering % center the table
    \caption{Pendulum swing-up optimal policy variables}. 
    \label{tb:optimalswingup}
    \begin{tabular}{p{1.0cm} p{2.5cm} p{1.0cm} p{1.5cm} }
        \hline \hline \noalign{\smallskip} \noalign{\smallskip} \noalign{\smallskip} \noalign{\smallskip}
        %%%%%%%%%%%%%%%%%%%%%%
        \textbf{Variable} & \textbf{Description} & \textbf{Units} & \textbf{Dimensions} \\ 
        %%%%%%%%%%%%%%%%%%%%%%
        \hline \hline \noalign{\smallskip} 
        \multicolumn{4}{c}{\textbf{Control inputs}}\\ \noalign{\smallskip}  \hline \hline
        \noalign{\smallskip} 
        %%%%%%%%%%%%%%%%%%%%%%
        $\tau$ & Actuator torque & $Nm$ & [$ML^2T^{-2}$]\\ 
        %%%%%%%%%%%%%%%%%%%%%%
        \hline \hline \noalign{\smallskip} 
        \multicolumn{4}{c}{\textbf{State variables}}\\ \noalign{\smallskip}  \hline \hline \noalign{\smallskip} 
        %%%%%%%%%%%%%%%%%%%%%%
        $\theta$ & Joint angle & $rad$ & []\\ \noalign{\smallskip} \hline \noalign{\smallskip}
        $\dot{\theta}$ & Joint angular velocity & $rad/sec$ & [$T^{-1}$] \\
        %%%%%%%%%%%%%%%%%%%%%%
        \hline \hline \noalign{\smallskip} 
        \multicolumn{4}{c}{\textbf{System parameters}}\\ \noalign{\smallskip}  \hline\hline  \noalign{\smallskip} 
        %%%%%%%%%%%%%%%%%%%%%%
        $m$ & Pendulum mass & $kg$ & [$M$]  \\ \noalign{\smallskip} \hline \noalign{\smallskip}
        $g$ & Gravity       & $m/s^2$ & [$LT^{-2}$]  \\ \noalign{\smallskip} \hline \noalign{\smallskip}
        $l$ & Pendulum lenght & $m$ & [$L$]  \\ \noalign{\smallskip} \hline \noalign{\smallskip}
        %%%%%%%%%%%%%%%%%%%%%%
        \hline \hline \noalign{\smallskip} 
        \multicolumn{4}{c}{\textbf{Problem parameters}}\\ \noalign{\smallskip}  \hline\hline  \noalign{\smallskip} 
        %%%%%%%%%%%%%%%%%%%%%%
        $q$ & Weight parameter  & $Nm$ & [$ML^2T^{-2}$]   \\ \noalign{\smallskip} \hline \noalign{\smallskip}
        $\tau_{max}$ & Maximum torque & $Nm$ & [$ML^2T^{-2}$] \\ \noalign{\smallskip} \hline \noalign{\smallskip}
        \hline \noalign{\smallskip}
        %\bottomrule[\heavyrulewidth] 
    \end{tabular}
\end{table}
%%%%%%%%%%%%%%%%%%%%%%%%%%%%%%%%%%%%%%%%%%%%
It is interesting to note that while there are three system parameters $m$, $g$, and $l$, they only appear independently in two groups in the dynamic equation. We can thus consider only two system parameters. For convenience, we selected $mgl$, corresponding to the maximum static gravitational torque (i.e., when the pendulum is horizontal) and the natural frequency $\omega=\sqrt{\frac{g}{l}}$, as listed in Table \ref{tb:2param}.
%%%%%%%%%%%%%%%%%%%%%%%%%%%%%%%%%%%%%%%%%%%%
\begin{table}[htb]
    \centering % center the table
    \caption{Pendulum reduced system parameters}. 
    \label{tb:2param}
    \begin{tabular}{p{1.0cm} p{2.5cm} p{1.0cm} p{1.5cm} }
        \hline \hline \noalign{\smallskip} \noalign{\smallskip} 
        \textbf{Variable} & \textbf{Description} & \textbf{Units} & \textbf{Dimensions} \\ \noalign{\smallskip}  \hline\hline  \noalign{\smallskip} 
        %%%%%%%%%%%%%%%%%%%%%%
        $mgl$ & Maximum gravitational torque  & $Nm$ & [$ML^2T^{-2}$]  \\ \noalign{\smallskip} \hline \noalign{\smallskip}
        $\omega = \sqrt{\frac{g}{l}}$ & Natural frequency & $sec^{-1}$ & [$T^{-1}$]  \\ \noalign{\smallskip} \hline \noalign{\smallskip}
        \hline \noalign{\smallskip}
        %\bottomrule[\heavyrulewidth] 
    \end{tabular}
\end{table}
%%%%%%%%%%%%%%%%%%%%%%%%%%%%%%%%%%%%%%%%%%%%

% %%%%%%%%%%%%%%%%%%%%%%%%%%%%%%%%%%%%%%%%%%%%
% \subsection{Dimensionless dynamics}

% %%%%%%%%%%%%%%%%%%%%%%
% \begin{align}
% %ml^2 \ddot{\theta} + mgl \sin \theta &= \tau  \\
% \frac{\ddot{\theta}}{\omega^2} + \sin \theta &= \frac{\tau}{mgl}
% \end{align}
% %%%%%%%%%%%%%%%%%%%%%%

% %%%%%%%%%%%%%%%%%%%%%%
% \begin{equation}
% m = (n = 5 ) - ( p = 2 ) = 3
% \end{equation}
% %%%%%%%%%%%%%%%%%%%%%%

% %%%%%%%%%%%%%%%%%%%%%%
% \begin{align}
% \Pi_1 &= \tau^* = \frac{\tau}{mgl} \quad \quad \frac{[ML^2T^{-2}]}{[M][LT^{-2}][L]} \\
% \Pi_2 &= \theta^* = \theta \quad \quad [-]\\
% \Pi_3 &= \ddot{\theta}^* = \frac{ \ddot{\theta}  }{ \omega^2 } \quad \quad \frac{[T^{-2}]}{[T^{-1}][T^{-1}]} 
% \end{align}
% %%%%%%%%%%%%%%%%%%%%%%

% %%%%%%%%%%%%%%%%%%%%%%
% \begin{equation}
% \tau^*
% =
% f\left(
% \theta^*,\ddot{\theta}^*
% \right) = \ddot{\theta}^* + \sin \theta^*
% \end{equation}
% %%%%%%%%%%%%%%%%%%%%%%

\subsubsection{Dimensional analysis}
\label{sec:dimanalpendulum}

Here, we have one control input, two states, two system parameters, and two task parameters, for a total of $1+(n=2)+(m=4)=7$ variables involved. In those variables, only $d=2$ independent dimensions ( $ML^2T^{-2}$ and $T^{-1}$ ) are present. Using $c_1 = mgl$ and $c_2 = \omega$ as the repeated variables leads to the following dimensionless groups:
%%%%%%%%%%%%%%%%%%%%%%
\begin{align}
\Pi_1 &= \tau^* = \frac{\tau}{mgl} \quad \quad \frac{[ML^2T^{-2}]}{[M][LT^{-2}][L]} \\
\Pi_2 &= \theta^* = \theta \quad \quad []\\
\Pi_3 &= \dot{\theta}^* = \frac{ \dot{\theta}  }{ \omega } \quad \quad \frac{[T^{-1}]}{[T^{-1}]} \\
\Pi_4 &= \tau_{max}^* = \frac{\tau_{max}}{mgl} \quad \quad \frac{[ML^2T^{-2}]}{[M][LT^{-2}][L]} \\
\Pi_5 &= q^* = \frac{q}{mgl} \quad \quad \frac{[ML^2T^{-2}]}{[M][LT^{-2}][L]} 
\end{align}
%%%%%%%%%%%%%%%%%%%%%%
All three torque variables ($\tau$, $q$, and $\tau_{max}$) are scaled by the maximum gravitational torque, and the pendulum velocity variable is scaled by the natural pendulum frequency. The transformation matrices are thus:
%%%%%%%%%%%%%%%%%%%%%%
\begin{align}
\tau^* &= 
\underbrace{\left[  1/mgl \right]}_{T_u}
\, \tau  \label{eq:Tupendulum} \\
\begin{bmatrix}
\theta^* \\ \dot{\theta}^*
\end{bmatrix} &= 
\underbrace{
\begin{bmatrix}
    1 & 0 \\ 0 & 1/\omega
\end{bmatrix}
}_{T_x} \, 
\begin{bmatrix}
\theta \\ \dot{\theta}
\end{bmatrix}
 \label{eq:Txpendulum} \\
\underbrace{
\begin{bmatrix}
q^* \\ \tau_{max}^*
\end{bmatrix} 
}_{c^*} 
&= 
\underbrace{
\begin{bmatrix}
 0 & 0  & 1/mgl & 0 \\  0 & 0  & 0 &  1/mgl
\end{bmatrix}
}_{T_c} \, 
\underbrace{
\begin{bmatrix}
mgl \\ \omega \\ q \\ \tau_{max}
\end{bmatrix}
}_{c} 
 \label{eq:Tcpendulum} 
\end{align}
%%%%%%%%%%%%%%%%%%%%%%
By applying the Buckingham $\pi$ theorem \cite{buckingham_physically_1914}, Equation \eqref{eq:pendulumpi} can be restated as a relationship between the five dimensionless $\Pi$ groups:
%%%%%%%%%%%%%%%%%%%%%%
\begin{equation}
\tau^*
=
\pi^* \left(
\theta, \dot{\theta}^*,
q^* , \tau_{max}^* 
\right)
\end{equation}
%%%%%%%%%%%%%%%%%%%%%%

\begin{figure}[ht]
\begin{center}
\includegraphics[width=0.95\linewidth]{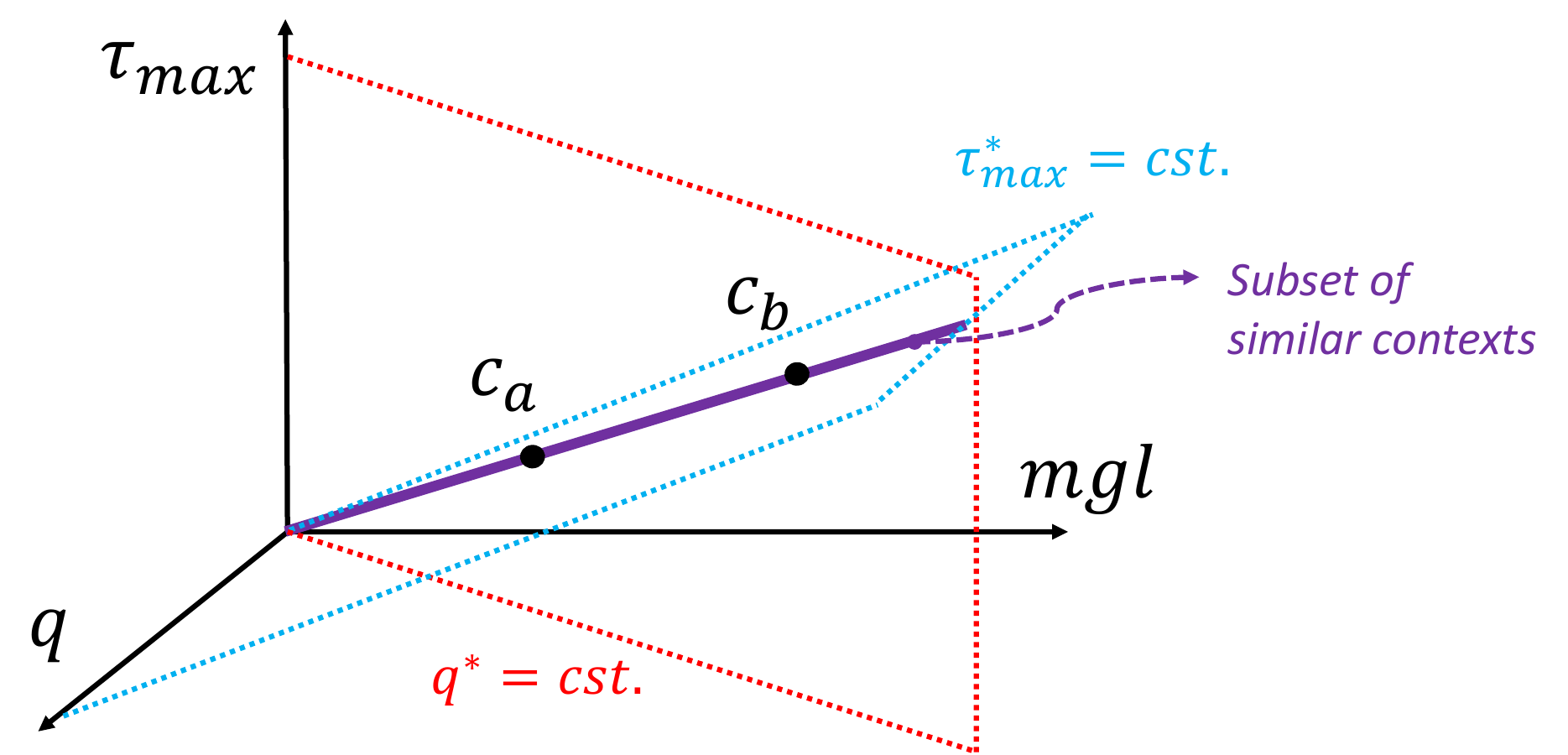}
\caption{The dimensionally similar subset (equal $c^*$) can be represented as a line in a 3D space for the pendulum swing-up problem. The feedback law solutions to problem with context variables on the same line are equivalent in dimensionless form.}
\label{fig:context_3d}
\end{center}
\end{figure}

According to the results of Section \ref{sec:dimenanalysis}, for dimensionally similar swing-up contexts (meaning those with equal $q^*$ and $\tau_{max}^*$ ratios), the optimal feedback laws should be equivalent in their dimensionless forms. In other words, the optimal policy $f_a$, found in the specific context $c_a = [m_a,l_a,g_a,q_a,\tau_{max,a}]$, and the optimal policy $f_b$, in a second context, $c_b = [m_b,l_b,g_b,q_b,\tau_{max,b}]$, are equal when restated in dimensionless form: $f_a^*=f_b^*$ if $q^*_a = q^*_b$ and $\tau_{max,a}^* = \tau_{max,b}^*$. Furthermore, $f_b$ can be obtained from $f_a$ or vice versa using the scaling formula given by Equation \eqref{eq:ab_transform} if this condition is met. However, if $q^*_a \neq q^*_b$ or $\tau_{max,a}^* \neq \tau_{max,b}^*$, then $f_a$ cannot provide us with information on $f_b$ without additional assumptions. Figure \ref{fig:context_3d} illustrates that for the pendulum swing-up problem the similarity condition can be represented as a line in a three dimension space created by three dimensional context variables. Each conditions of equal values of $q^*$ and $\tau_{max}^*$, is a plane in this space, and the intersection of the two plane is the subset of context meeting the two conditions. Also, it is interesting to note that the fourth context variable $\omega$ is not an additional axis here because it is not involved in eq. \eqref{eq:Tcpendulum}.

%%%%%%%%%%%%%%%%%%%%%%
%%%%%%%%%%%%%%%%%%%%%%

%%%%%%%%%%%%%%%%%%%%%%
\subsubsection{Numerical results}
Here, we use a numerical algorithm (methodological details are presented in Section \ref{sec:metho}) to compute numerical solutions to the motion control problem defined by Equations \eqref{eq:pendulum_dynamics}, \eqref{eq:pendulum_cost}, and \eqref{eq:pendulum_constraints}. The algorithm computes feedback laws in the form of look-up tables, based on a discretized grid of the state space. The optimal (up to discretization errors) feedback laws are computed for nine instances of context variables, which are listed in Table \ref{tb:9contexts}. In those nine contexts, there are three subsets of three dimensionally similar contexts. Also, each subset includes the same three pendulums: a regular pendulum, one that is two times longer, and one that is twice as heavy (as illustrated in Figure \ref{fig:big_picture}). Contexts $c_a$, $c_b$, and $c_c$ describe a task where the torque is limited to half the maximum gravitational torque. Contexts $c_d$, $c_e$, and $c_f$ describe a task where the application of large torques is highly penalized by the cost function. Contexts $c_g$, $c_h$, and $c_i$ describe a task where position errors are highly penalized by the cost function.   
%%%%%%%%%%%%%%%%%%%%%%%%%%%%%%%%%%%%%%%%%%%%
\begin{table}[h]
    \centering % center the table
    \caption{Pendulum swing-up problem context variables.} 
    \label{tb:9contexts}
    \begin{tabular}{ p{2.0cm} p{0.8cm} p{0.8cm} p{0.8cm} p{0.8cm} p{0.8cm} }
        \hline \hline \noalign{\smallskip} \noalign{\smallskip} 
        %%%%%%%%%%%%%%%%%%%%%%
        & $m$ & $g$ & $l$ & $q$ & $\tau_{max}$ \\ \hline
        %%%%%%%%%%%%%%%%%%%%%
        %%%%%%%%%%%%%%%%%%%%%%
        \hline \hline \noalign{\smallskip} 
        \multicolumn{6}{c}{\textbf{Problems with $\tau_{max}^* = 0.5$ and $q^* = 0.1$} }\\ \noalign{\smallskip}  \hline\hline  \noalign{\smallskip} 
        %%%%%%%%%%%%%%%%%%%%%%
        Context $c_a$ : & 1.0 & 10.0 & 1.0 & 1.0 & 5.0 \\
        Context $c_b$ : & 1.0 & 10.0 & 2.0 & 2.0 & 10.0 \\
        Context $c_c$ : & 2.0 & 10.0 & 1.0 & 2.0 & 10.0 \\
        %%%%%%%%%%%%%%%%%%%%
        \hline \hline \noalign{\smallskip} 
        \multicolumn{6}{c}{\textbf{Problems with $\tau_{max}^* = 1.0$ and $q^* = 0.05$} }\\ \noalign{\smallskip}  \hline\hline  \noalign{\smallskip} 
        %%%%%%%%%%%%%%%%%%%%%%
        Context $c_d$ : & 1.0 & 10.0 & 1.0 & 0.5 & 10.0 \\
        Context $c_e$ : & 1.0 & 10.0 & 2.0 & 1.0 & 20.0 \\
        Context $c_f$ : & 2.0 & 10.0 & 1.0 & 1.0 & 20.0 \\
        %%%%%%%%%%%%%%%%%%%%%
        \hline \hline \noalign{\smallskip} 
        \multicolumn{6}{c}{\textbf{Problems with $\tau_{max}^* = 1.0$ and $q^* = 10$} }\\ \noalign{\smallskip}  \hline\hline  \noalign{\smallskip} 
        %%%%%%%%%%%%%%%%%%%%%%
        Context $c_g$ : & 1.0 & 10.0 & 1.0 & 100.0 & 10.0 \\
        Context $c_h$ : & 1.0 & 10.0 & 2.0 & 200.0 & 20.0 \\
        Context $c_i$ : & 2.0 & 10.0 & 1.0 & 200.0 & 20.0 \\
        %%%%%%%%%%%%%%%%%%%%%
        \hline \hline
    \end{tabular}
\end{table}
%%%%%%%%%%%%%%%%%%%%%%%%%%%%%%%%%%%%%%%%%%%%

\begin{figure*}[htp]
        \centering
        \vspace{-10pt}
        \subfloat[Feedback law $f$]{\includegraphics[width=6cm]{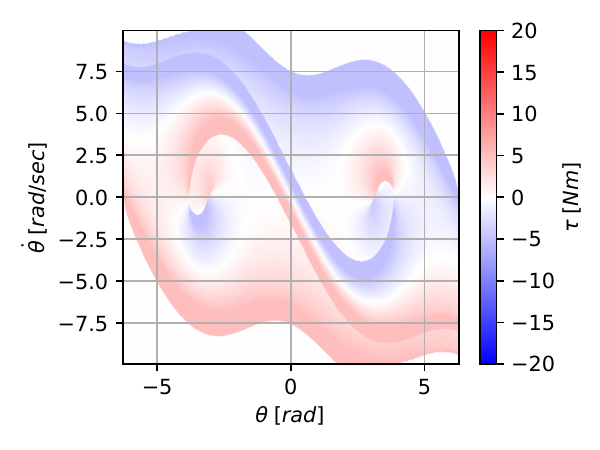}}
        \subfloat[Dimensionless feedback law $f^*$]{\includegraphics[width=6cm]{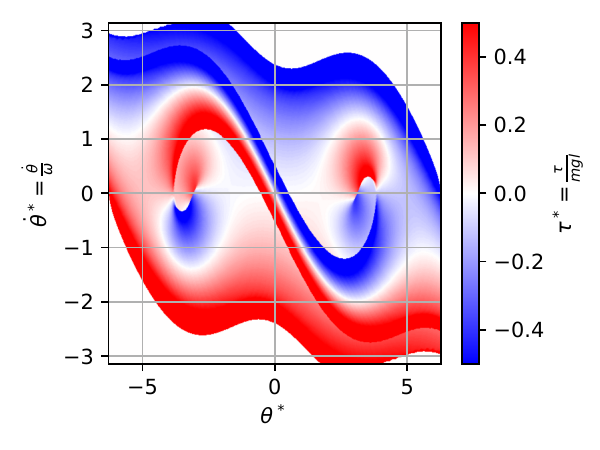}}   
        \subfloat[Optimal trajectory, starting at $\theta=-\pi$]{\includegraphics[width=6cm]{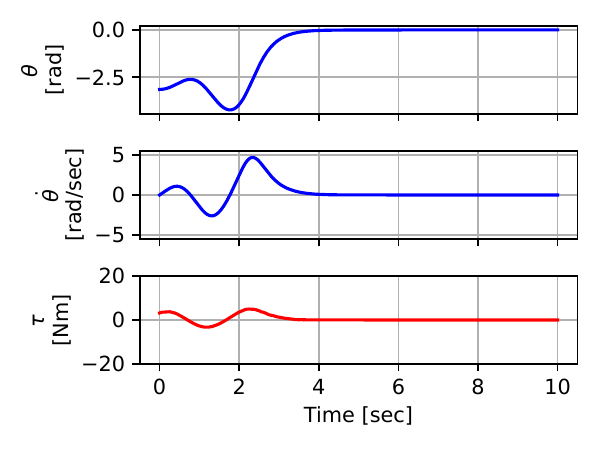}}
    
    \caption{Numerical results for context $c_a$}.
    \label{fig:c1}
\end{figure*}

\begin{figure*}[htp]
        \centering
        \vspace{-10pt}
        \subfloat[Feedback law $f$]{\includegraphics[width=6cm]{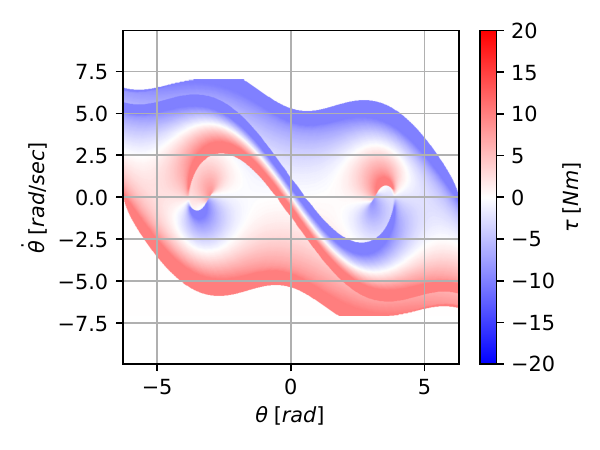}}
        \subfloat[Dimensionless feedback law $f^*$]{\includegraphics[width=6cm]{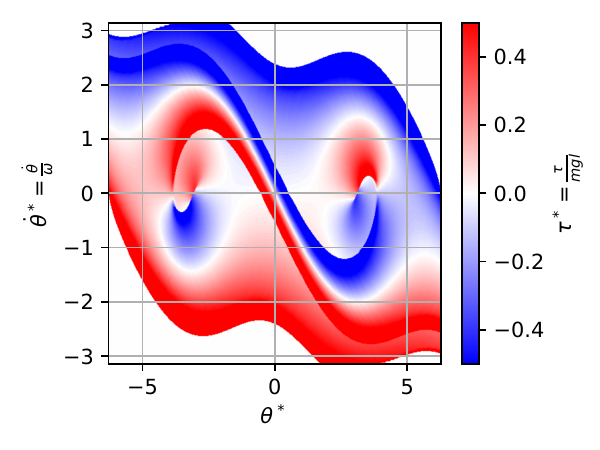}}   
        \subfloat[Optimal trajectory, starting at $\theta=-\pi$]{\includegraphics[width=6cm]{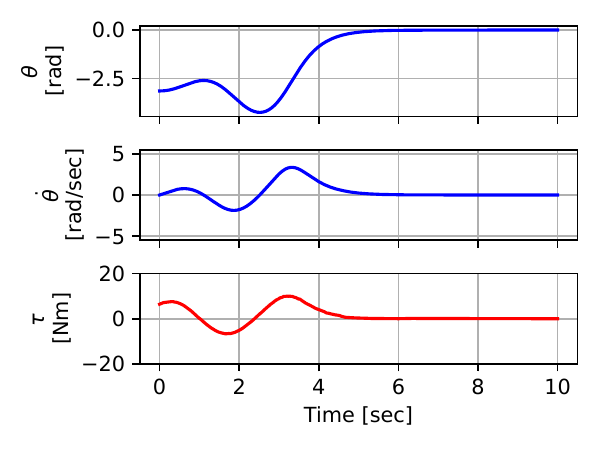}}
    
    \caption{Numerical results for context $c_b$}.
    \label{fig:c2}
\end{figure*}

\begin{figure*}[htp]
    
        \centering
        \vspace{-10pt}
        \subfloat[Feedback law $f$]{\includegraphics[width=6cm]{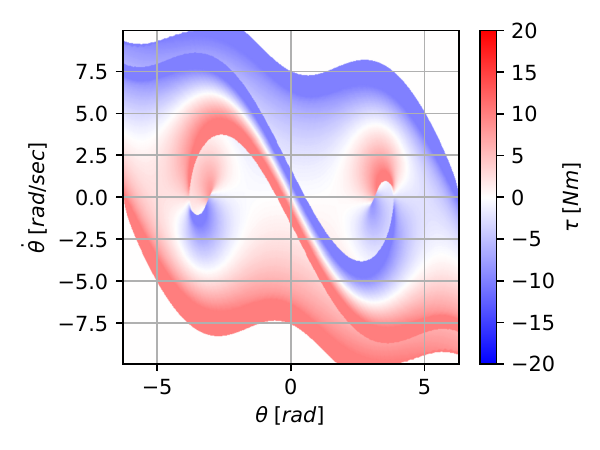}}
        \subfloat[Dimensionless feedback law $f^*$]{\includegraphics[width=6cm]{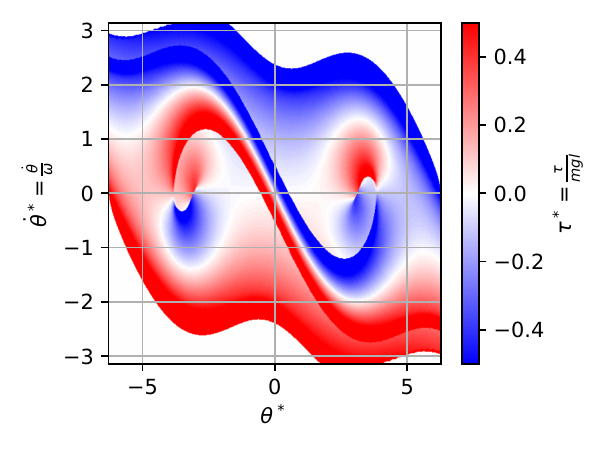}}   
        \subfloat[Optimal trajectory, starting at $\theta=-\pi$]{\includegraphics[width=6cm]{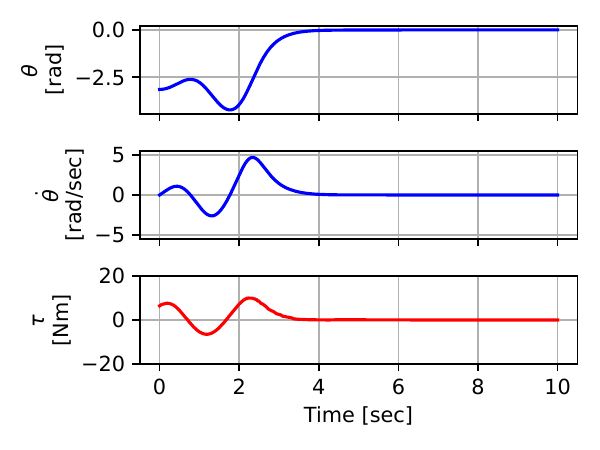}}
    
    \caption{Numerical results for context $c_c$.}
    \label{fig:c3}
\end{figure*}

\begin{figure*}[htp]
    
        \centering
        \vspace{-10pt}
        \subfloat[Feedback law $f$]{\includegraphics[width=6cm]{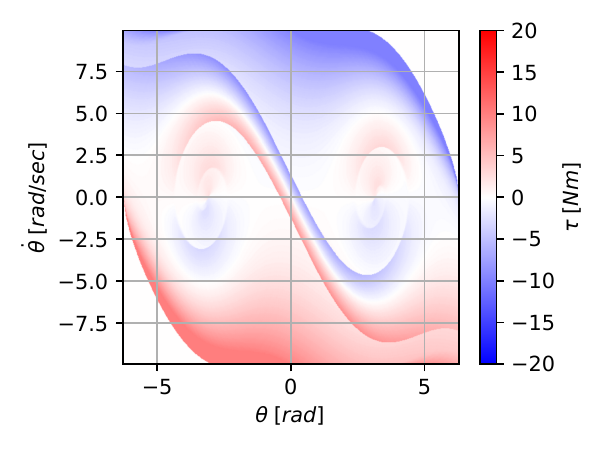}}
        \subfloat[Dimensionless feedback law $f^*$]{\includegraphics[width=6cm]{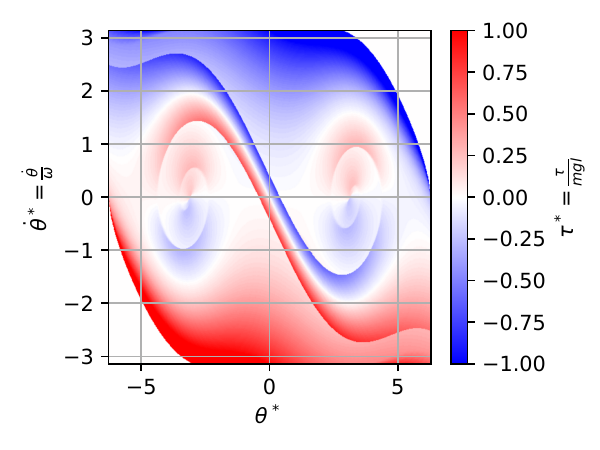}}   
        \subfloat[Optimal trajectory, starting at $\theta=-\pi$]{\includegraphics[width=6cm]{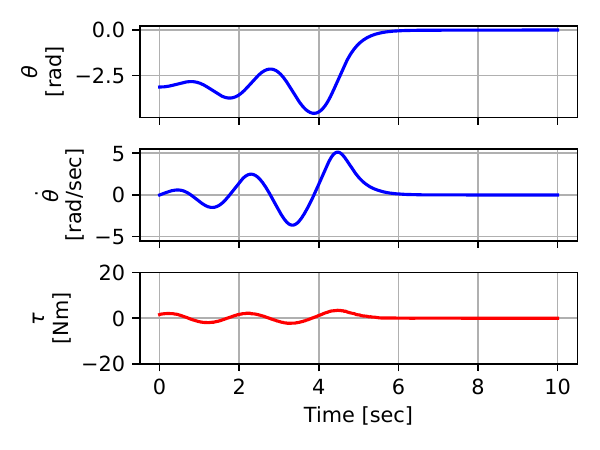}}
    
    \caption{Numerical results for context $c_d$}.
    \label{fig:c4}
\end{figure*}

\begin{figure*}[htp]
    
        \centering
        \vspace{-10pt}
        \subfloat[Feedback law $f$]{\includegraphics[width=6cm]{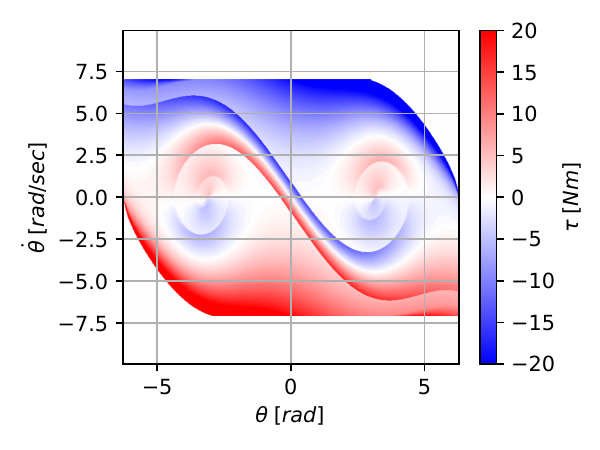}}
        \subfloat[Dimensionless feedback law $f^*$]{\includegraphics[width=6cm]{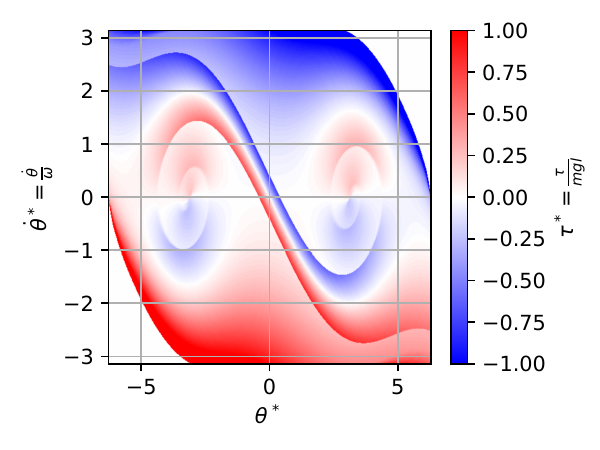}}   
        \subfloat[Optimal trajectory, starting at $\theta=-\pi$]{\includegraphics[width=6cm]{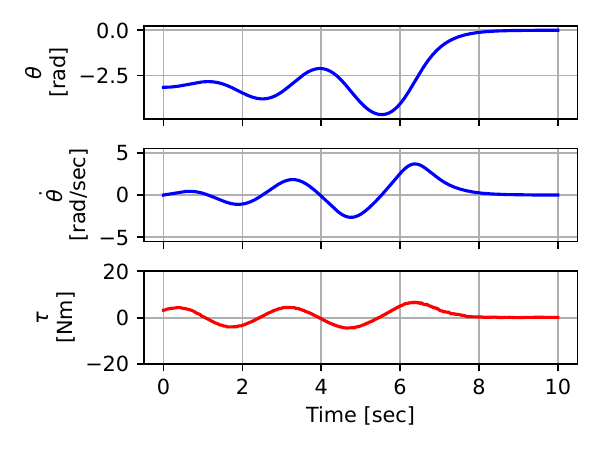}}
    
    \caption{Numerical results for context $c_e$}.
    \label{fig:c5}
\end{figure*}

\begin{figure*}[htp]
    
        \centering
        \vspace{-10pt}
        \subfloat[Feedback law $f$]{\includegraphics[width=6cm]{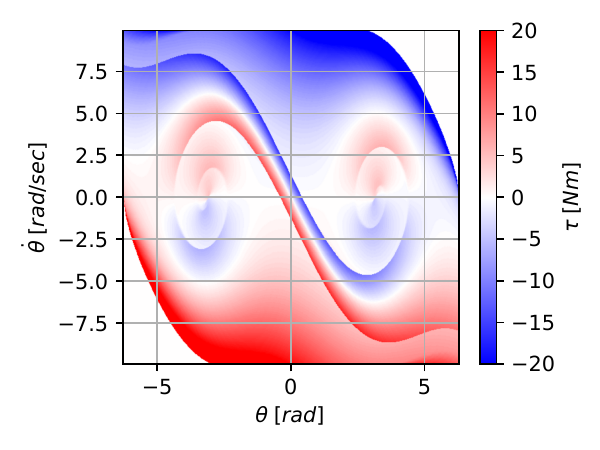}}
        \subfloat[Dimensionless feedback law $f^*$]{\includegraphics[width=6cm]{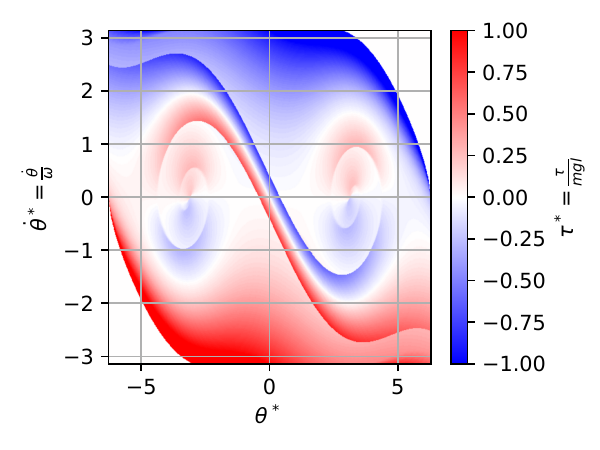}}   
        \subfloat[Optimal trajectory, starting at $\theta=-\pi$]{\includegraphics[width=6cm]{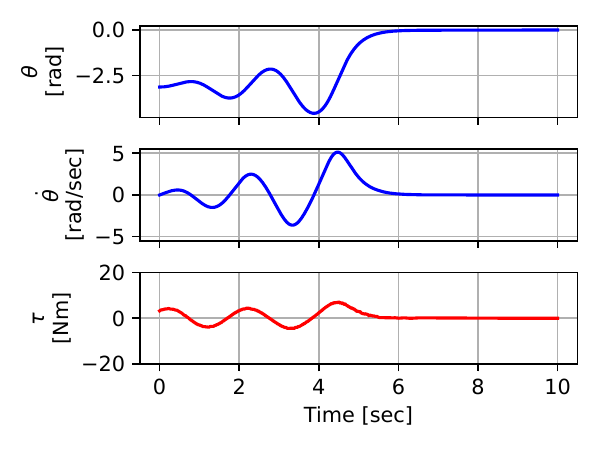}}
    
    \caption{Numerical results for context $c_f$}.
    \label{fig:c6}
\end{figure*}

\begin{figure*}[htp]
    
        \centering
        \vspace{-10pt}
        \subfloat[Feedback law $f$]{\includegraphics[width=6cm]{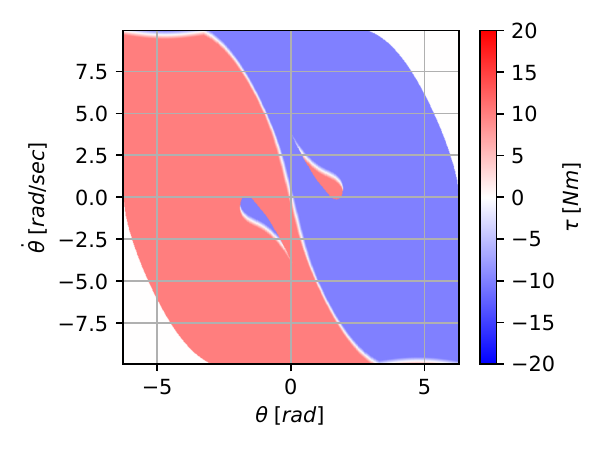}}
        \subfloat[Dimensionless feedback law $f^*$]{\includegraphics[width=6cm]{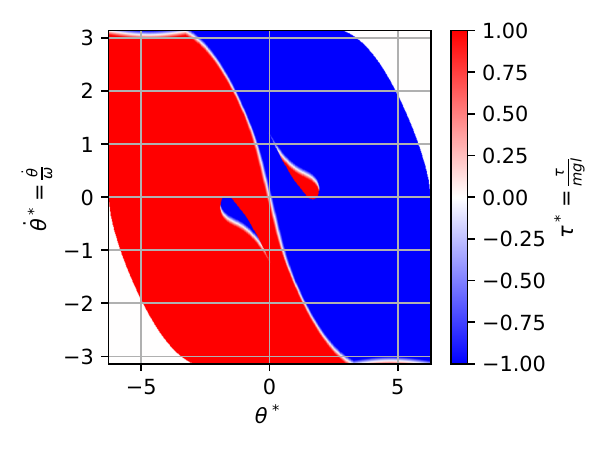}}   
        \subfloat[Optimal trajectory, starting at $\theta=-\pi$]{\includegraphics[width=6cm]{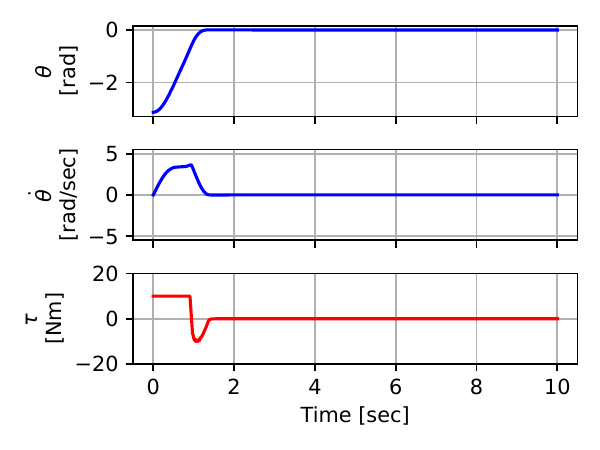}}
    
    \caption{Numerical results for context $c_g$}.
    \label{fig:c7}
\end{figure*}

\begin{figure*}[htp]
    
        \centering
        \vspace{-10pt}
        \subfloat[Feedback law $f$]{\includegraphics[width=6cm]{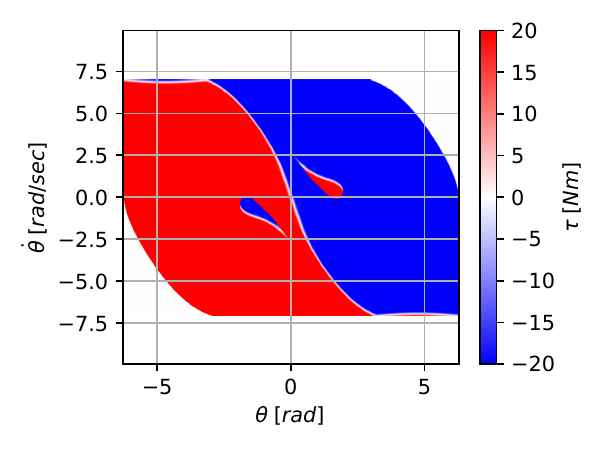}}
        \subfloat[Dimensionless feedback law $f^*$]{\includegraphics[width=6cm]{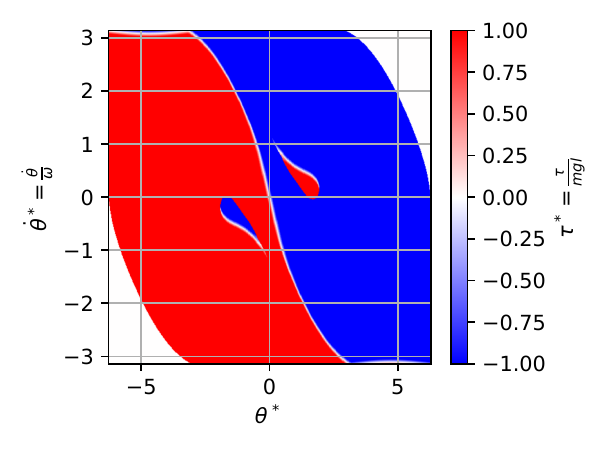}}   
        \subfloat[Optimal trajectory, starting at $\theta=-\pi$]{\includegraphics[width=6cm]{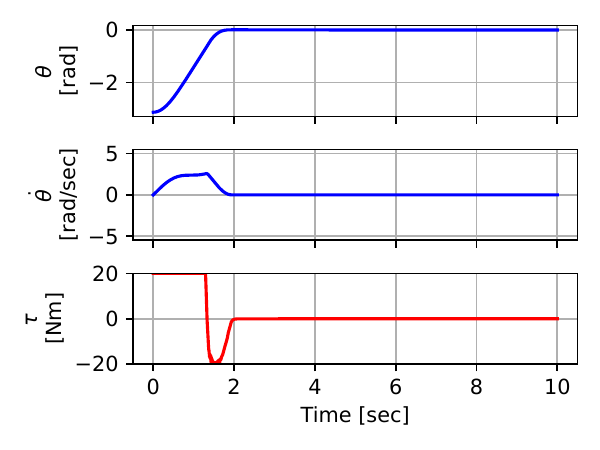}}
    
    \caption{Numerical results for context $c_h$}.
    \label{fig:c8}
\end{figure*}

\begin{figure*}[tp]
    
        \centering
        \vspace{-10pt}
        \subfloat[Feedback law $f$]{\includegraphics[width=6cm]{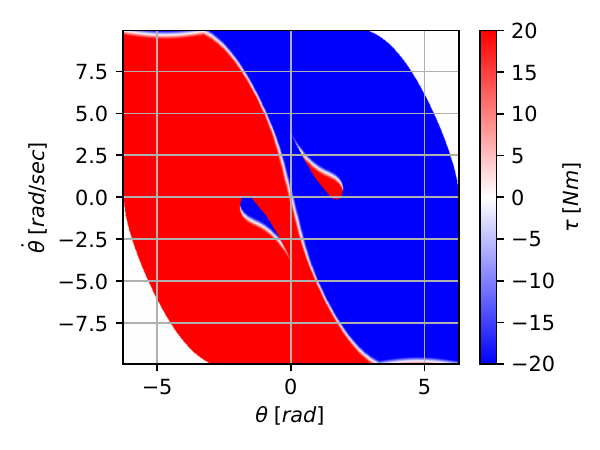}}
        \subfloat[Dimensionless feedback law $f^*$]{\includegraphics[width=6cm]{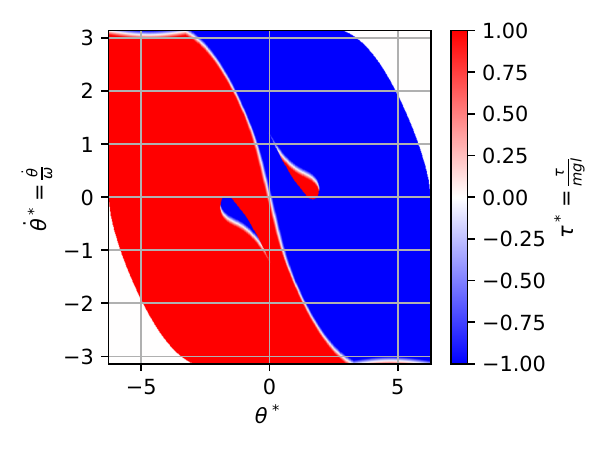}}   
        \subfloat[Optimal trajectory, starting at $\theta=-\pi$]{\includegraphics[width=6cm]{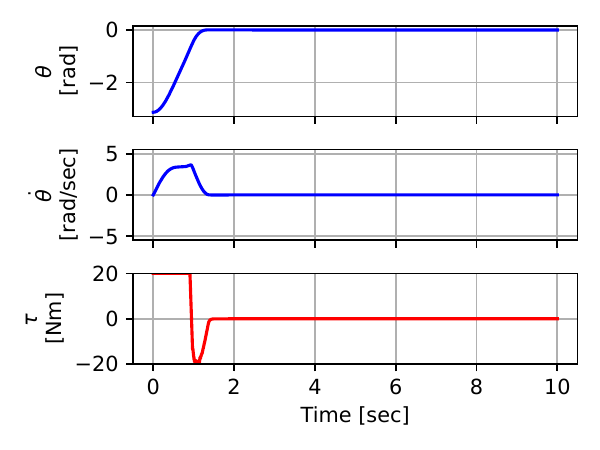}}
    
    \caption{Numerical results for context $c_i$}.
    \label{fig:c9}
\end{figure*}

Figures \ref{fig:c1} to \ref{fig:c9} illustrate that, for each subset with equal dimensionless context, the dimensional feedback laws generated look numerically very similar. They are similar up to the scaling of their axis, if we neglect slight differences due to discretization errors. Furthermore, the figures also illustrate that the dimensionless version of the feedback laws ($f^*$), computed using Equation \eqref{eq:f2star}, are equal within each dimensionally similar subset. These were the expected results predicted by the dimensional analysis presented in Section \ref{sec:dimenanalysis}.

%%%%%%%%%%%%%%%%%%%%%%%%%%%%
%%%%%%%%%%%%%%%%%%%%%%%%%%%%
%%%%%%%%%%%%%%%%%%%%%%%%%%%%
%%%%%%%%%%%%%%%%%%%%%%%%%%%%
%%%%%%%%%%%%%%%%%%%%%%%%%%%%
%%%%%%%%%%%%%%%%%%%%%%%%%%%
%%%%%%%%%%%%%%%%%%%%%%%%%%%%
%%%%%%%%%%%%%%%%%%%%%%%%%%%%
%%%%%%%%%%%%%%%%%%%%%%%%%%%%
%%%%%%%%%%%%%%%%%%%%%%%%%%%%
In terms of how this can be applied in a practical scenario, we see that if we compute the feedback law given in Figure \ref{fig:c1}(a), we can obtain the feedback law given in Figure \ref{fig:c2}(a) directly by scaling the original policy with Equation \eqref{eq:ab_transform}, using the appropriate context variables, without having to recompute. In some sense, Equation \eqref{eq:ab_transform} provides us with the ability to adjust the feedback law spontaneously to conform with new system parameters $mgl$ or $\omega$, as would be the case with an analytical solution, even when working with black box results in the form of a table look-up. But the equivalence of the scaled solution is only guaranteed within a dimensionally similar context subset, which is the main limitation of this approach. The feedback law given in Figure \ref{fig:c1}(a) cannot be scaled into the feedback law given in Figure \ref{fig:c7}(a), for instance, since $\tau^*_{max}$ and $q^*$ are not equals. It is also interesting to note that trajectory solutions from the same starting point—and computed cost-to-go functions (not illustrated)—are also all equivalent, up to scaling factors, within similar subgroups. Hence, optimal trajectories and cost-to-go solutions could also be shared and transferred between similar systems using the same technique that we demonstrate here for feedback laws.

\subsubsection{Regimes of solutions}

\begin{figure}[htp]
    \begin{center}
        \includegraphics[width=0.9\linewidth]{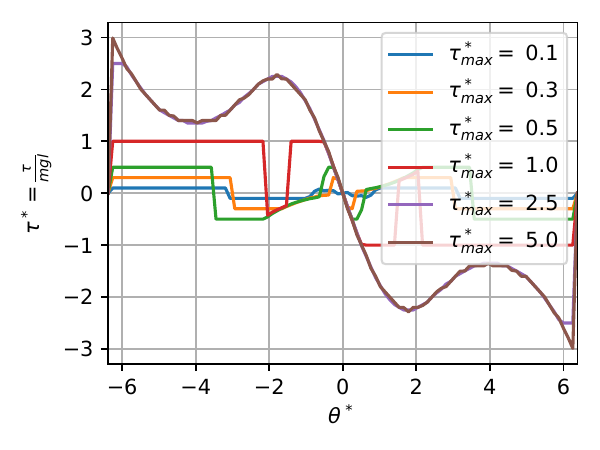}
        \caption{Optimal dimensionless policy for various contexts: $\tau^* = \pi^*( \theta^* , \dot{\theta}^* = 0 , q^* = 0.5 , \tau^*_{max} = [0.1, ... , 5.0] )$.}
        \label{fig:torque_sensitivity}
    \end{center}
\end{figure}

\begin{figure}[htp]
    \begin{center}
        \includegraphics[width=0.9\linewidth]{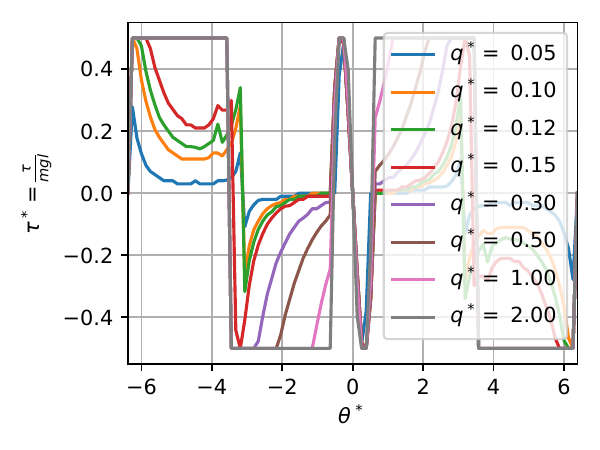}
        \caption{Optimal dimensionless policy for various contexts:  $\tau^* = \pi^*( \theta^*  , \dot{\theta}^* = 0 , q^* = [0.05, ... , 2.0]  , \tau^*_{max} = 0.5 ).$}
        \label{fig:q_sensitivity}
    \end{center}
\end{figure}

In some situations, changing a context variable will not have any effect on the optimal policy. For instance, for a torque-limited optimal pendulum swing-up problem, augmenting $\tau_{max}$ or $q$ while keeping the other value fixed will have little effect above a given threshold. For instance, if we look at the solutions for contexts $c_d$, $c_e$, and $c_f$, using a large amount of torque is so highly penalized by the cost function that the saturation limit does not have much impact on the solution (except for edge cases on the boundary). Thus, we would expect that augmenting $\tau_{max}$ would not change the solution. Figures \ref{fig:torque_sensitivity} and \ref{fig:q_sensitivity} show a slice (to allow for visualization) of the dimensionless optimal policy solution for various contexts. Figure \ref{fig:torque_sensitivity} illustrates the results of changing $\tau_{max}^*$ while keeping $q^*$ fixed. We can see that when $\tau_{max}^*<0.3$, the policy is almost always on the min--max allowable torque values; this behavior is often called \textit{bang--bang}. At the other extreme, when 
$\tau_{max}^*>2.5$, the policy solution is continuous and almost never affected by the saturation. Figure \ref{fig:q_sensitivity} illustrates the results of changing $q^*$ while keeping $\tau_{max}^*$ fixed. We can see that when $q^*<0.1$, the optimal policy solution does not reach min--max saturation, while when $q^*>1.0$, the policy is almost always on the min--max allowable values.
%%%%%%%%%%%%%%%%%%%%%%
%%%%%%%%%%%%%%%%%%%%%%
%%%%%%%%%%%%%%%%%%%%%%
%%%%%%%%%%%%%%%%%%%%%%
%%%%%%%%%%%%%%%%%%%%%%
\begin{figure}[ht]
    \begin{center}
        \includegraphics[width=0.9\linewidth]{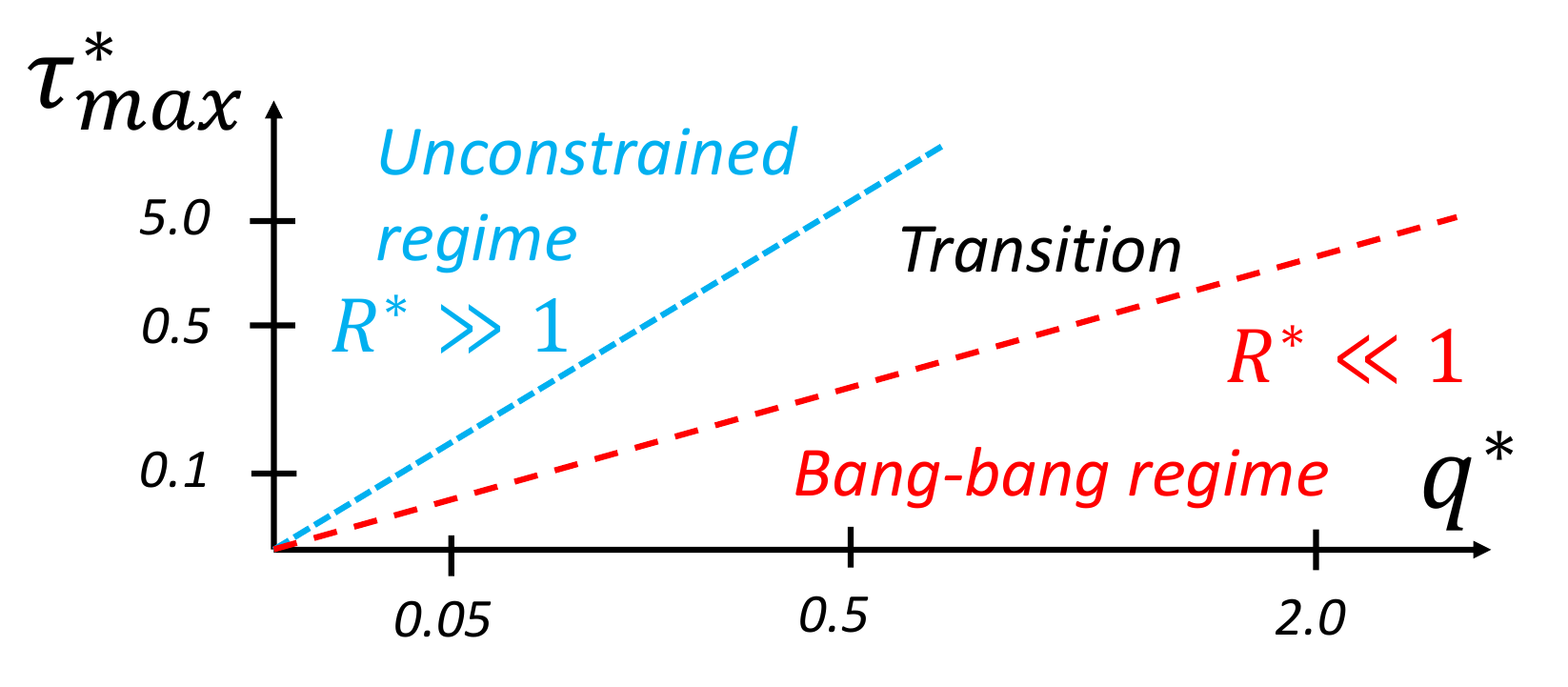}
        \caption{Regime zones for a torque-limited pendulum swing-up problem.}
        \label{fig:regimeszones}
    \end{center}
\end{figure}
%%%%%%%%%%%%%%%%%%%%%%

We can see that, for extreme context values, two types of behavior occur, illustrated as regions in the dimensionless context space in Figure \ref{fig:regimeszones}. Those regions are best characterized by a ratio of $q^*$ and $\tau_{max}^*$, a new dimensionless value that we define as the ratio of the maximum torque saturation $\tau_{max}$ over the weight parameter in the cost function $q$:
%%%%%%%%%%%%%%%%%%%%%%
\begin{align}
R^* = \frac{\tau^*_{max}}{q^*} = \frac{\tau_{max}}{q} \quad\quad 
\end{align}
%%%%%%%%%%%%%%%%%%%%%%
When the value of $R^* \approx 1$, the policy solution is partially continuous and reaches the min--max value in some other region of the state space, this is a behavior we call the transition regime. When the value of $R^* \ll 1$, the constraint on torque drives the solution to exhibit bang--bang behavior. In this region (that we approximate here, based on our sensitivity analysis, as $R^* \leq 0.1$), the global policy is only a function of $\tau_{max}^*$:
%%%%%%%%%%%%%%%%%%%%%%
\begin{align}
\pi^*( \theta^* , \dot{\theta}^*, q^*, \tau_{max}^*) &\approx 
\pi^*( \theta^* , \dot{\theta}^* , \tau_{max}^* ) \; \text{if} \; R^* \ll 1
\label{eq:bangbang_policy}
\end{align}
%%%%%%%%%%%%%%%%%%%%%%
i.e., the value of $q^*$ does not affect the solution. On the other hand, when the value of $R^* \gg 1$, the policy is unconstrained. In this region (that we approximate here, based on our sensitivity analysis, as $R^* \geq 10$), the global policy is only a function of $q^*$ since the constraint is so far away:
%%%%%%%%%%%%%%%%%%%%%%
\begin{align}
\pi^*( \theta^* , \dot{\theta}^*, q^*, \tau_{max}^*) &\approx 
\pi^*( \theta^* , \dot{\theta}^* , q^* ) \; \text{if} \;  R^* \gg 1 
\label{eq:unconstrained_policy}
\end{align}
%%%%%%%%%%%%%%%%%%%%%%
The concept of regime is often leveraged in fluid mechanics. It allows us to generalize results between situations where the relevant dimensionless numbers do not match exactly. For instance, when the Mach number is small ($Ma < 0.3$), we can generally assume there to be in an incompressible regime where various speeds of sound would not change the behavior much. Here, for the purpose of transferring policy solutions between contexts, this means that the condition of having the same exact dimensionless context variables can be relaxed with an inequality that corresponds to a regime. For instance, if we have two contexts in the unconstrained regime, it is sufficient to match only $q^*$ to create equivalent dimensionless policies.
%if we consider Equation \eqref{eq:unconstrained_policy}, the condition of having equivalent dimensionless feedback laws is relaxed with an inequality for one of the context variables, as follows:
\begin{Proposition}
If it is assumed that Equation \eqref{eq:unconstrained_policy} holds, the condition of having equivalent dimensionless feedback laws is relaxed to an inequality for one of the context variables, as follows:
%%%%%%%%%%%%%%%%%%%%%%
\begin{align}
&f_a^*( \theta^* , \dot{\theta}^* ) \approx f_b^*( \theta^* , \dot{\theta}^*) 
\\ &\text{if} \quad q^*_a=q^*_b \quad \text{\textit{and}} \quad R^*_a \gg 1 \quad \text{\textit{and}} \quad R^*_b \gg 1
\end{align}
%%%%%%%%%%%%%%%%%%%%%%
\end{Proposition}
\begin{proof}
First, if $R_a^* \gg 1$ and $R_b^* \gg 1$ then from Equation \eqref{eq:unconstrained_policy} we can approximate the policy not to be a function of $\tau_{max}^*$:
%%%%%%%%%%%%%%%%%%%%%%
\begin{align}
f_a^*(\theta^* , \dot{\theta}^*) &= \pi^*( \theta^* , \dot{\theta}^*, q_a^*, \tau_{max,a}^*) \approx 
\pi^*( \theta^* , \dot{\theta}^* , q_a^* ) \\
f_b^*(\theta^* , \dot{\theta}^*) &= \pi^*( \theta^* , \dot{\theta}^*, q_b^*, \tau_{max,b}^*) \approx 
\pi^*( \theta^* , \dot{\theta}^* , q_b^* ) 
\end{align}
%%%%%%%%%%%%%%%%%%%%%%
Hence, if $q_a^*=q_b^*$ we have:
%%%%%%%%%%%%%%%%%%%%%%
\begin{align}
f_a^*(\theta^* , \dot{\theta}^*)  \approx 
\pi^*( \theta^* , \dot{\theta}^* , q_a^* ) = \pi^*( \theta^* , \dot{\theta}^* , q_b^* ) \approx 
f_b^*(\theta^* , \dot{\theta}^*) 
\end{align}
%%%%%%%%%%%%%%%%%%%%%%
\end{proof}
Also, for two contexts in a bang--bang regime, it is sufficient to match only $\tau_{max}^*$ to have equivalent dimensionless policies.
\begin{Proposition}
If it is assumed that Equation \eqref{eq:bangbang_policy} holds, the condition of having equivalent dimensionless feedback laws is relaxed to an inequality for one of the context variables, as follows:
%%%%%%%%%%%%%%%%%%%%%%
\begin{align}
&f_a^*( \theta^* , \dot{\theta}^* ) \approx f_b^*( \theta^* , \dot{\theta}^*) 
\\ &\text{if} \quad \tau^*_{max,a}=\tau^*_{max,b} \quad \text{\textit{and}} \quad R^*_a \ll 1 \quad \text{\textit{and}} \quad R^*_b \ll 1
\end{align}
%%%%%%%%%%%%%%%%%%%%%%
\end{Proposition}
\begin{proof}
First, if $R_a^* \ll 1$ and $R_b^* \ll 1$ then from Equation \eqref{eq:bangbang_policy} we can approximate the policy not to be a function of $q^*$:
%%%%%%%%%%%%%%%%%%%%%%
\begin{align}
f_a^*(\theta^* , \dot{\theta}^*) &= \pi^*( \theta^* , \dot{\theta}^*, q_a^*, \tau_{max,a}^*) \approx 
\pi^*( \theta^* , \dot{\theta}^* , \tau^*_{max,a} ) \\
f_b^*(\theta^* , \dot{\theta}^*) &= \pi^*( \theta^* , \dot{\theta}^*, q_b^*, \tau_{max,b}^*) \approx 
\pi^*( \theta^* , \dot{\theta}^* , \tau^*_{max,b} ) 
\end{align}
%%%%%%%%%%%%%%%%%%%%%%
Hence, if $\tau^*_{max,a}=\tau^*_{max,a}$ we have:
%%%%%%%%%%%%%%%%%%%%%%
\begin{align}
f_a^*(\theta^* , \dot{\theta}^*)  &\approx 
\pi^*( \theta^* , \dot{\theta}^* , \tau^*_{max,a} ) \\
f_a^*(\theta^* , \dot{\theta}^*)  &\approx  \pi^*( \theta^* , \dot{\theta}^* , \tau^*_{max,b} ) \\
f_a^*(\theta^* , \dot{\theta}^*)  &\approx f_b^*(\theta^* , \dot{\theta}^*) 
\end{align}
%%%%%%%%%%%%%%%%%%%%%%
\end{proof}

From another point of view, assuming that one of those regimes applies means that we could have removed one variable from the context at the start of the dimensional analysis. All in all, the impact of identifying such regimes is that we can increase the size of the context subset to which the dimensionless version of the policy should be equivalent, leading to a potentially larger pool of systems that can share a learned policy and numerical results.

\subsubsection{Methodology}
\label{sec:metho}
We obtained the optimal feedback law presented in this section using a basic dynamic programming algorithm \cite{bertsekas_dynamic_2012} on a discretized version of the continuous system. The approach is almost equivalent to the value iteration  algorithm \cite{sutton_reinforcement_2018}—which is sometimes referred to as model-based reinforcement learning—with the exception that, here, the total number of iteration steps was fixed (corresponding to a very long time horizon approximating an infinite horizon), instead of the iteration being stopped after reaching a convergence criterion. This approach was chosen to enable the collection of consistent results across all contexts that lead to a wide range of order-of-magnitude cost-to-go solutions. The time step was set to 0.025 s, the state space was discretized into an even 501 x 501 grid, and the continuous torque input was discretized into 101 discrete control options. Special out-of-bounds and on-target termination states were included to guarantee convergence \cite{bertsekas_dynamic_2012}. Also, using dynamic programming made the setting of additional parameters to define the domain necessary. Although those parameters should not affect the optimal policy far away from the boundaries, dimensionless versions of those parameters were kept fixed in all the experiments, as follows:
%%%%%%%%%%%%%%%%%%%%%%
\begin{align}
\theta^*_{max} &= \theta_{max} = 2 \pi \\
\dot{\theta}^*_{max} &= \frac{ \dot{\theta}_{max} }{\omega} = \pi \\
t^*_{f} &= t_{f} \; \omega = 20 \times 2 \pi 
\end{align}
%%%%%%%%%%%%%%%%%%%%%%
where $\theta_{max}$ is the range of angles for which the optimal policy is solved, set to one full revolution; $\dot{\theta}_{max}$ is the range of angular velocity for which the optimal policy is solved; and $t_{f}$ is the time horizon, set to 20 periods of the pendulum using the natural frequency. The source code is available online at the following link: \url{https://github.com/alx87grd/DimensionlessPolicies}, and this Google Colab page allows users to reproduce the results: \url{https://colab.research.google.com/drive/1kf3apyHlf5t7XzJ3uVM8mgDsneVK_63r?usp=sharing}.

%%%%%%%%%%%%%%%%%%%%%%
\subsection{Optimal motion for a longitudinal car on a slippery surface}
\label{sec:optimalcar}
The second numerical example is a simplified car positioning task. We use this example to illustrate that an optimal feedback law in the form of a table look-up generated for a car of a given size, can be transferred to a car of a different size if the motion control problem is dimensionally similar. The example includes state constraints and a different type of non-linearity (i.e. is its not similar to the pendulum swing-up) to illustrate how generic the developed dimensionless polices concept are.

%%%%%%%%%%%%%%%%%%%%%%
\begin{figure}[ht]
    \begin{center}
        \includegraphics[width=0.9\linewidth]{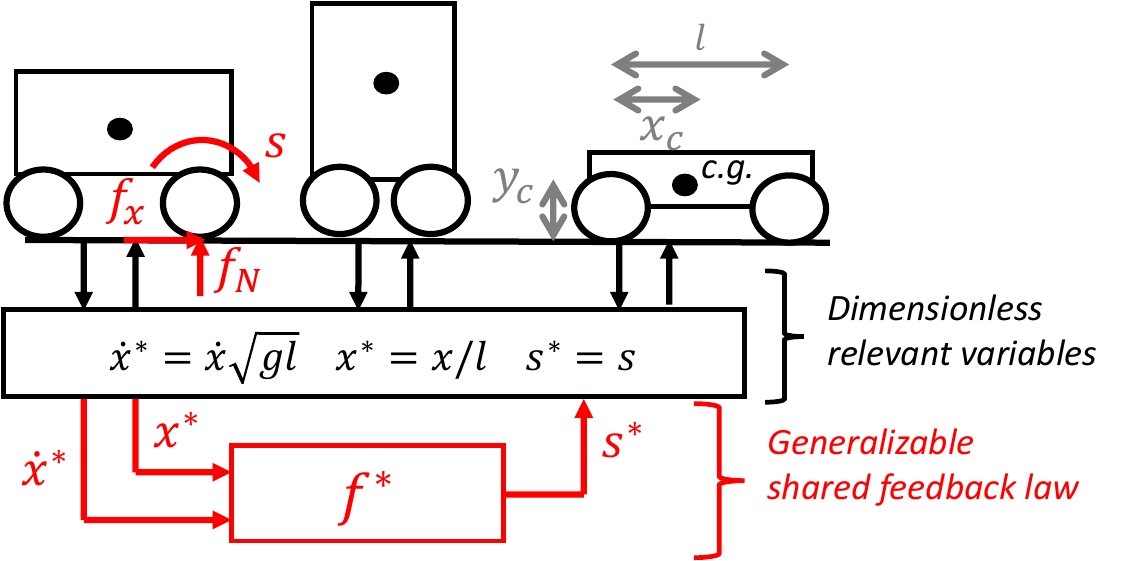}
        \caption{Car positioning motion control problem.}
        \label{fig:cars}
    \end{center}
\end{figure}
%%%%%%%%%%%%%%%%%%%%%%

\subsubsection{Motion control problem}
The motion control problem is defined here as finding a feedback law to control the dynamic system, as described by the following differential equation:
%%%%%%%%%%%%%%%%%%%%%%
\begin{equation}
\ddot{x} = \frac{ \mu(s) g x_c }{ l + \mu(s) y_c } 
\quad \quad \text{with} \quad \quad 
\mu(s) = \frac{ f_x }{ f_n^{f} }  =  \frac{ 2 }{ 1 + e^{-70s} } - 1
\label{eq:car_dynamics}
\end{equation}
%%%%%%%%%%%%%%%%%%%%%%
where $\mu(s)$ is the ratio of vertical to horizontal forces on the front wheel, that is, a non-linear function of the front wheel slip $s$. The above equations represent a simple dynamic model of the longitudinal motion of a car, assuming that the controller can impose the wheel slip of the front wheel and that suspensions are infinitely rigid (but that weight transfer is included). Interestingly, it is already standard practice to model the ground--tire interaction with an empirical curve $\mu(s)$ relating two dimensionless variables. 

The objective is to minimize the infinite horizon quadratic cost function given by:
% %%%%%%%%%%%%%%%%%%%%%%
% \begin{equation}
% J = \int{\left( q^{-2} x^2 + 0 \, \dot{x}^2 + 1 \, s^2 \right) dt }
% \label{eq:car_cost}
% \end{equation}
% %%%%%%%%%%%%%%%%%%%%%%
%%%%%%%%%%%%%%%%%%%%%
\begin{equation}
J = \int_0^{\infty}{ \left( q^{-2} x^2 + \, s^2 \right) dt }
\label{eq:car_cost}
\end{equation}
%%%%%%%%%%%%%%%%%%%%%%
subject to the constraints of keeping ground reaction forces positive, as given by:
%%%%%%%%%%%%%%%%%%%%%%
\begin{align}
0 \leq& f_n^{f} = g x_c - \ddot{x} y_c
\\
0 \leq& f_n^{r} = g (l-x_c) + \ddot{x} y_c
\label{eq:car_constraints}
\end{align}
%%%%%%%%%%%%%%%%%%%%%%
where the weight transfer potentially limits the allowable motions. Note that the cost function parameter $q$ in this problem has a power of minus two to have a value with units of length, and all parameters are time-independent constants. The solution to this problem, i.e., the optimal policy for all contexts, involves the variables listed in Table \ref{tb:optimalcar} and should be of the form given by:
%%%%%%%%%%%%%%%%%%%%%%
\begin{equation}
\underbrace{s}_{\text{input}}
=
\pi \left(
\underbrace{ x, \dot{x} }_{\text{states}},
\underbrace{ x_c , y_c , g , l , q 
}_{\text{Context $c$}}
\right)
\label{eq:carpi}
\end{equation}
%%%%%%%%%%%%%%%%%%%%%
.
%%%%%%%%%%%%%%%%%%%%%%%%%%%%%%%%%%%%%%%%%%%%
\begin{table}[htb]
    \centering % center the table
    \caption{Longitudinal car optimal policy variables}. 
    \label{tb:optimalcar}
    \begin{tabular}{p{1.0cm} p{2.5cm} p{1.0cm} p{1.5cm}}
        \hline \hline \noalign{\smallskip} \noalign{\smallskip} \noalign{\smallskip} \noalign{\smallskip}
        %%%%%%%%%%%%%%%%%%%%%%
        \textbf{Variable} & \textbf{Description} & \textbf{Units} & \textbf{Dimensions} \\ 
        %%%%%%%%%%%%%%%%%%%%%%
        \hline \hline \noalign{\smallskip} 
        \multicolumn{4}{c}{\textbf{Control inputs}}\\ \noalign{\smallskip}  \hline \hline
        \noalign{\smallskip} 
        %%%%%%%%%%%%%%%%%%%%%%
        $s$ & Wheel slip & - & []\\ 
        %%%%%%%%%%%%%%%%%%%%%%
        \hline \hline \noalign{\smallskip} 
        \multicolumn{4}{c}{\textbf{State variables}}\\ \noalign{\smallskip}  \hline \hline \noalign{\smallskip} 
        %%%%%%%%%%%%%%%%%%%%%%
        $x$ & Car position & $m$ & [L]\\ \noalign{\smallskip} \hline \noalign{\smallskip}
        $\dot{x}$ & Car velocity & $m/sec$ & [$LT^{-1}$] \\
        %%%%%%%%%%%%%%%%%%%%%%
        \hline \hline \noalign{\smallskip} 
        \multicolumn{4}{c}{\textbf{System parameters}}\\ \noalign{\smallskip}  \hline\hline  \noalign{\smallskip} 
        %%%%%%%%%%%%%%%%%%%%%%
        $g$ & Gravity       & $m/s^2$ & [$LT^{-2}$]  \\ \noalign{\smallskip} \hline \noalign{\smallskip}
        $l$ & Length (wheel base) & $m$ & [$L$]  \\ \noalign{\smallskip} \hline \noalign{\smallskip}
        $x_c$ & center of gravity (CG) horizontal position& $m$ & [$L$]  \\ \noalign{\smallskip} \hline \noalign{\smallskip}
        $y_c$ & center of gravity (CG) vertical position & $m$ & [$L$]  \\ \noalign{\smallskip} \hline \noalign{\smallskip}
        %%%%%%%%%%%%%%%%%%%%%%
        \hline \hline \noalign{\smallskip} 
        \multicolumn{4}{c}{\textbf{Problem parameters}}\\ \noalign{\smallskip}  \hline\hline  \noalign{\smallskip} 
        %%%%%%%%%%%%%%%%%%%%%%
        $q$ & Weight parameter  & $m$ & [$L$]   \\ \noalign{\smallskip} \hline \noalign{\smallskip}
        %\bottomrule[\heavyrulewidth] 
    \end{tabular}
\end{table}
%%%%%%%%%%%%%%%%%%%%%%%%%%%%%%%%%%%%%%%%%%%%

\subsubsection{Dimensional analysis}

Here, we have one control input, two states, and five context parameters, for a total of $1+(n=2)+(m=5)=8$ variables. Of those variables, only $d=2$ independent dimensions (length $[L]$ and time $[T]$) are present. Using $c_1 = g$ and $c_2 = l$ as the repeated variables leads to the following dimensionless groups:
%%%%%%%%%%%%%%%%%%%%%%
\begin{align}
\Pi_1 &= s^* = s \quad \quad [] \\
\Pi_2 &= x^* = \frac{x}{l}  \quad \quad \frac{[L]}{[L]}\\
\Pi_3 &= \dot{x}^* = \frac{ \dot{x}  }{ \sqrt{gl} } \quad \quad \frac{[LT^{-1}]}{[LT^{-2}]^{1/2}[L]^{1/2}} \\
\Pi_4 &= x_c^* = \frac{x_c}{l}  \quad \quad \frac{[L]}{[L]}\\
\Pi_5 &= y_c^* = \frac{y_c}{l}  \quad \quad \frac{[L]}{[L]}\\
\Pi_6 &= q^* = \frac{q}{l}  \quad \quad \frac{[L]}{[L]}
\end{align}
%%%%%%%%%%%%%%%%%%%%%%
All three length variables are scaled by the wheel base, and the velocity variable is scaled using a combination of the wheel base and gravity. The transformation matrices are then as follows:
%%%%%%%%%%%%%%%%%%%%%%
\begin{align}
s* &= 
\underbrace{\left[  1 \right]}_{T_u}
\, s  \label{eq:Tupcar} \\
\begin{bmatrix}
x^* \\ \dot{x}^*
\end{bmatrix} &= 
\underbrace{
\begin{bmatrix}
    \frac{1}{l} & 0 \\ 0 & \frac{1}{\sqrt{gl}}
\end{bmatrix}
}_{T_x} \, 
\begin{bmatrix}
x \\ \dot{x}
\end{bmatrix}
 \label{eq:Txcar} \\
\underbrace{
\begin{bmatrix}
x_c^* \\ y_c^* \\ q^* 
\end{bmatrix} 
}_{c^*} 
&= 
\underbrace{
\begin{bmatrix}
 0 & 0  & 1/l & 0 & 0 \\
 0 & 0  & 0 &  1/l & 0 \\
  0 & 0  & 0 &  0 & 1/l
\end{bmatrix}
}_{T_c} \, 
\underbrace{
\begin{bmatrix}
g \\ l \\ x_c \\ y_c \\ q
\end{bmatrix}
}_{c} 
 \label{eq:Tccar} 
\end{align}
%%%%%%%%%%%%%%%%%%%%%%
By applying the Buckingham $\pi$ theorem \cite{buckingham_physically_1914}, Equation \eqref{eq:carpi} can be restated as a relationship between the six dimensionless $\Pi$ groups, as follows:
%%%%%%%%%%%%%%%%%%%%%%
\begin{equation}
s^*
=
\pi^* \left(
x^*, \dot{x}^*,
x_c^* , y_c^* , q^* 
\right)
\end{equation}
%%%%%%%%%%%%%%%%%%%%%%

\subsubsection{Numerical results}

Here, as in the pendulum example, numerical solutions to the motion control problem are computed for the nine instances of context variables listed in Table \ref{tb:9contextscar}. In those nine contexts, there are three subsets of three dimensionally similar contexts. Contexts $c_a$, $c_b$, and $c_c$ describe situations where the CG. horizontal position is at half the wheel base; contexts $c_d$, $c_e$ and $c_f$ describe situations in which the the CG is very high (and hence the cars are very limited by the weight transfer); and contexts $c_h$, $c_i$, and $c_j$ describe situations in which position errors are highly penalized by the cost function plus cars with a very low CG relative to the wheel base.
%%%%%%%%%%%%%%%%%%%%%%%%%%%%%%%%%%%%%%%%%%%%
\begin{table}[h]
    \centering % center the table
    \caption{Car problem parameters.} 
    \label{tb:9contextscar}
    \begin{tabular}{ p{2.0cm} p{0.8cm} p{0.8cm} p{0.8cm} p{0.8cm} p{0.8cm} }
        \hline \hline \noalign{\smallskip} \noalign{\smallskip} 
        %%%%%%%%%%%%%%%%%%%%%%
        & $l$ & $g$ & $x_c$ & $y_c$ & $q$ \\ \hline
        %%%%%%%%%%%%%%%%%%%%%
        %%%%%%%%%%%%%%%%%%%%%%
        \hline \hline \noalign{\smallskip} 
        \multicolumn{6}{c}{\textbf{Problems with $x_c^* = 0.5$, $y_c^* = 0.5$, and $q^* = 20$} }\\ \noalign{\smallskip}  \hline\hline  \noalign{\smallskip} 
        %%%%%%%%%%%%%%%%%%%%%%
        Context $c_a$ : & 2.0 & 9.8 & 1.0 & 1.0 & 40 \\
        Context $c_b$ : & 1.0 & 9.8 & 0.5 & 0.5 & 20 \\
        Context $c_c$ : & 3.0 & 9.8 & 1.5 & 1.5 & 60 \\
        %%%%%%%%%%%%%%%%%%%%
        \hline \hline \noalign{\smallskip} 
        \multicolumn{6}{c}{\textbf{Problems with $x_c^* = 0.5$, $y_c^* = 1.5$, and $q^* = 10$} }\\ \noalign{\smallskip}  \hline\hline  \noalign{\smallskip} 
        %%%%%%%%%%%%%%%%%%%%%%
        Context $c_d$ : & 2.0 & 9.8 & 1.0 & 3.0 & 20 \\
        Context $c_e$ : & 1.0 & 9.8 & 0.5 & 1.5 & 10 \\
        Context $c_f$ : & 3.0 & 9.8 & 1.5 & 4.5 & 30 \\
        %%%%%%%%%%%%%%%%%%%%%
        \hline \hline \noalign{\smallskip} 
        \multicolumn{6}{c}{\textbf{Problems with $x_c^* = 0.5$, $y_c^* = 0.1$, and $q^* = 2$ }}\\ \noalign{\smallskip}  \hline\hline  \noalign{\smallskip} 
        %%%%%%%%%%%%%%%%%%%%%%
        Context $c_g$ : & 2.0 & 9.8 & 1.0 & 0.2 & 4 \\
        Context $c_h$ : & 1.0 & 9.8 & 0.5 & 0.1 & 2 \\
        Context $c_i$ : & 3.0 & 9.8 & 1.5 & 0.3 & 6 \\
        %%%%%%%%%%%%%%%%%%%%%
        \hline \hline
    \end{tabular}
\end{table}
%%%%%%%%%%%%%%%%%%%%%%%%%%%%%%%%%%%%%%%%%%%%

\begin{figure*}[htp]
			
				\centering
				\vspace{-10pt}
				\subfloat[Feedback law $f$]{\includegraphics[width=6cm]{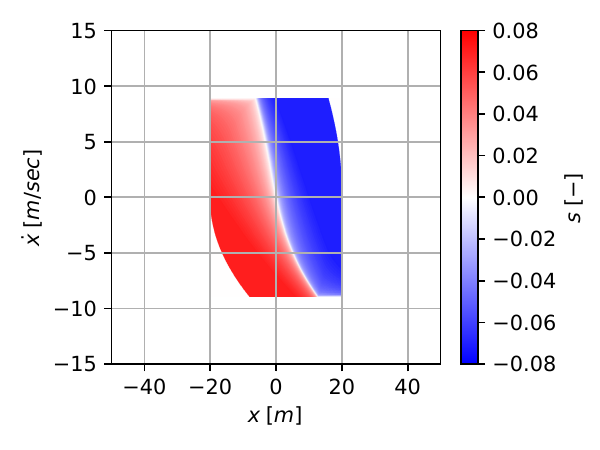}}
				\subfloat[Dimensionless feedback law $f^*$]{\includegraphics[width=6cm]{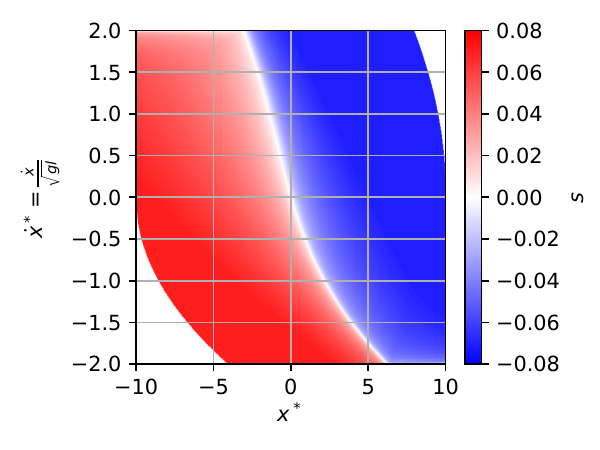}}   
				\subfloat[Optimal trajectory, starting at $x = -5l$]{\includegraphics[width=6cm]{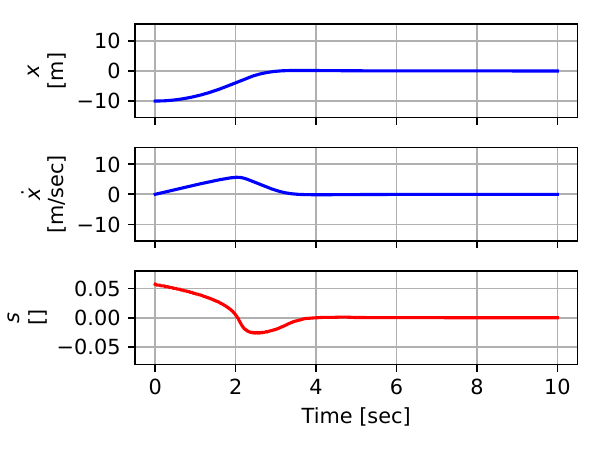}}
			
			\caption{Numerical results for context $c_a$.}
			\label{fig:c1c}
		\end{figure*}

\begin{figure*}[htp]
			
				\centering
				\vspace{-10pt}
				\subfloat[Feedback law $f$]{\includegraphics[width=6cm]{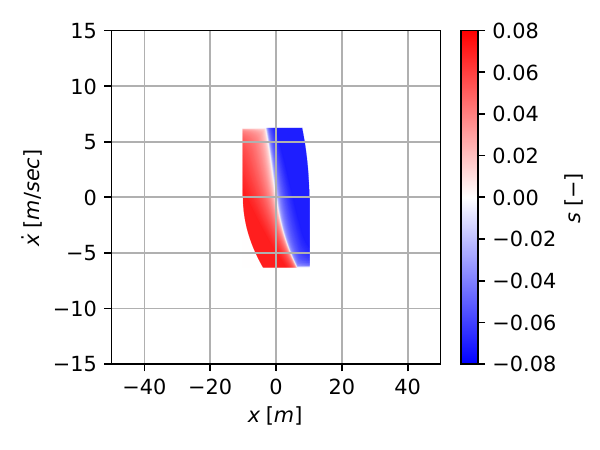}}
				\subfloat[Dimensionless feedback law $f^*$]{\includegraphics[width=6cm]{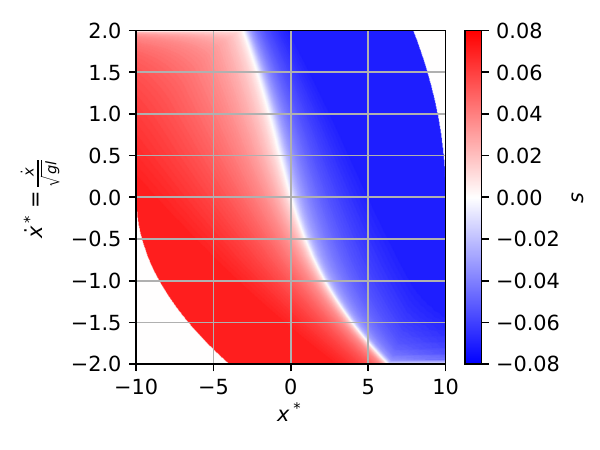}}   
				\subfloat[Optimal trajectory, starting at $x = -5 l$]{\includegraphics[width=6cm]{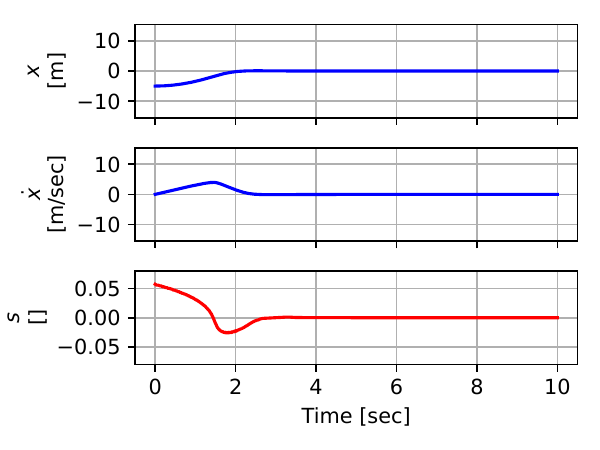}}
			
			\caption{Numerical results for context $c_b$.}
			\label{fig:c2c}
		\end{figure*}

\begin{figure*}[htp]
			
				\centering
				\vspace{-10pt}
				\subfloat[Feedback law $f$]{\includegraphics[width=6cm]{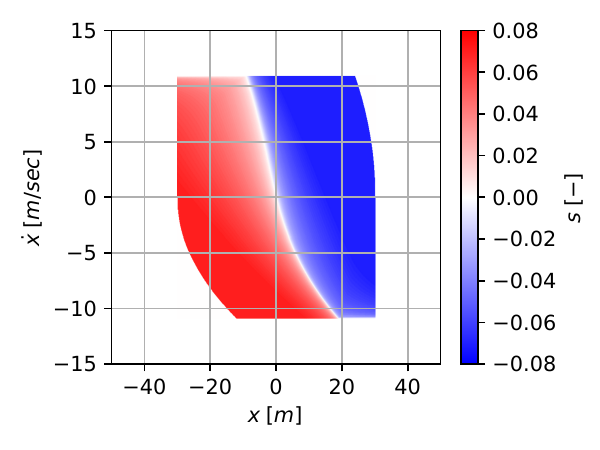}}
				\subfloat[Dimensionless feedback law $f^*$]{\includegraphics[width=6cm]{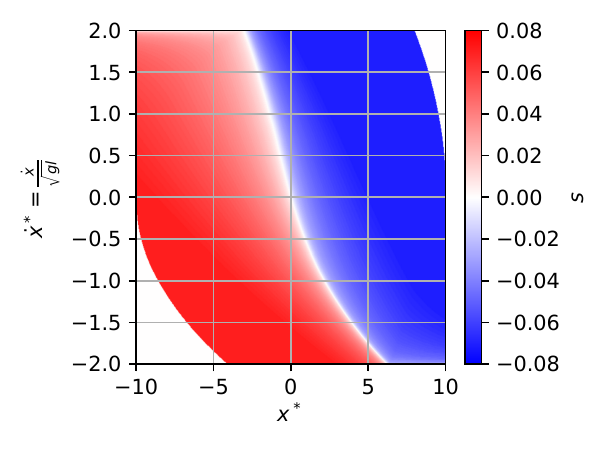}}   
				\subfloat[Optimal trajectory, starting at $x = -5 l$]{\includegraphics[width=6cm]{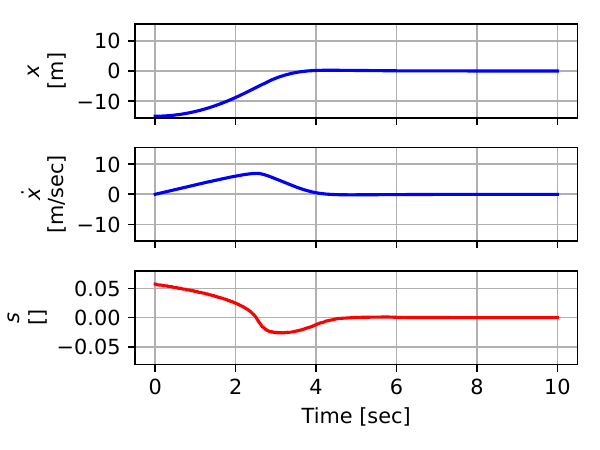}}
			
			\caption{Numerical results for context $c_c$.}
			\label{fig:c3c}
		\end{figure*}

\begin{figure*}[htp]
			
				\centering
				\vspace{-10pt}
				\subfloat[Feedback law $f$]{\includegraphics[width=6cm]{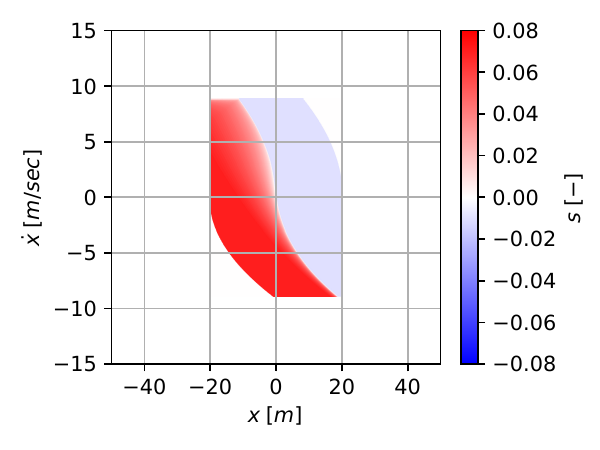}}
				\subfloat[Dimensionless feedback law $f^*$]{\includegraphics[width=6cm]{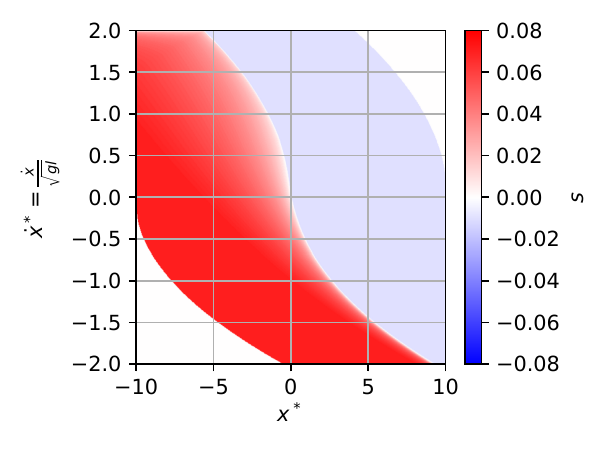}}   
				\subfloat[Optimal trajectory, starting at $x = -5 l$]{\includegraphics[width=6cm]{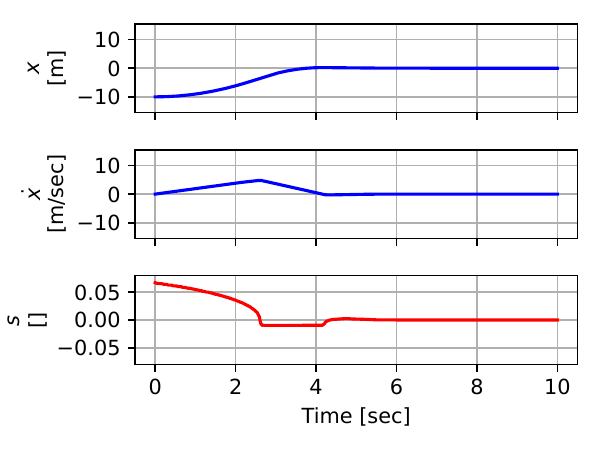}}
			
			\caption{Numerical results for context $c_d$.}
			\label{fig:c4c}
		\end{figure*}

\begin{figure*}[htp]
			
				\centering
				\vspace{-10pt}
				\subfloat[Feedback law $f$]{\includegraphics[width=6cm]{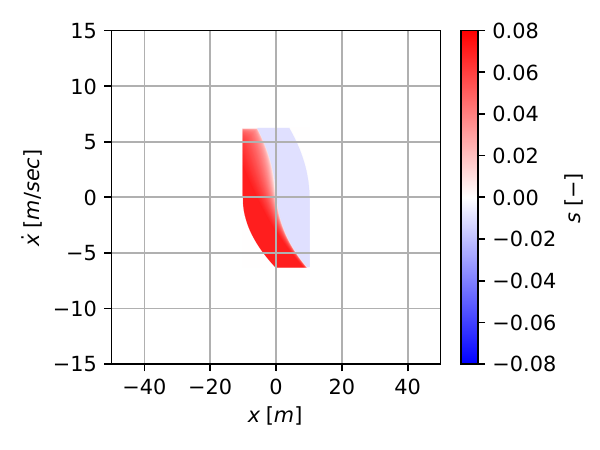}}
				\subfloat[Dimensionless feedback law $f^*$]{\includegraphics[width=6cm]{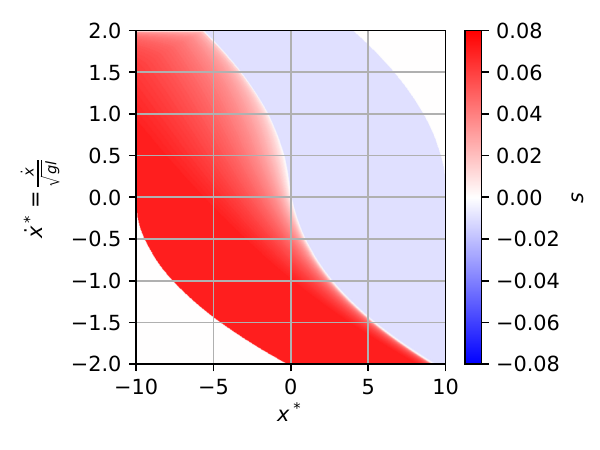}}   
				\subfloat[Optimal trajectory, starting at $x = -5 l$]{\includegraphics[width=6cm]{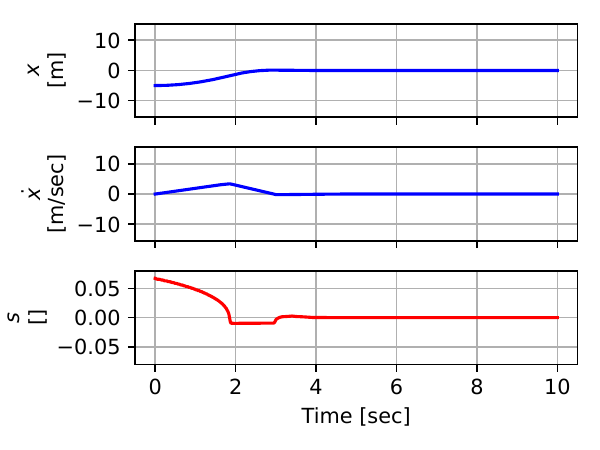}}
			
			\caption{Numerical results for context $c_e$.}
			\label{fig:c5c}
		\end{figure*}

\begin{figure*}[htp]
			
				\centering
				\vspace{-10pt}
				\subfloat[Feedback law $f$]{\includegraphics[width=6cm]{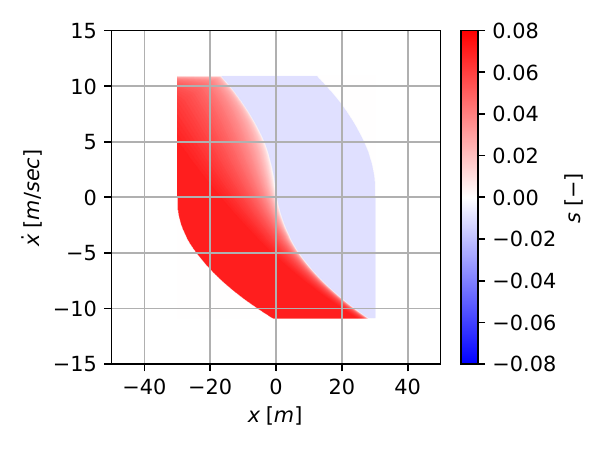}}
				\subfloat[Dimensionless feedback law $f^*$]{\includegraphics[width=6cm]{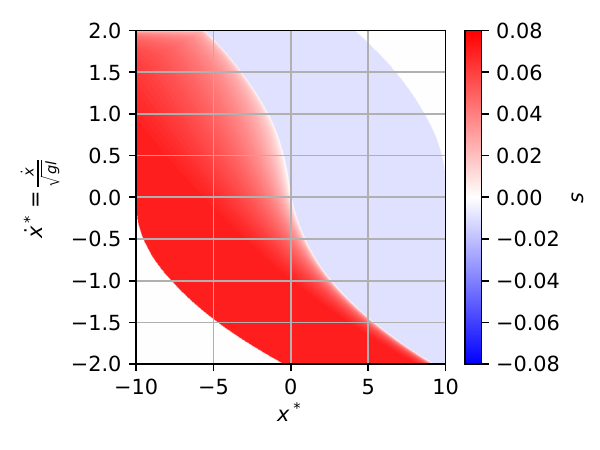}}   
				\subfloat[Optimal trajectory, starting at $x = -5 l$]{\includegraphics[width=6cm]{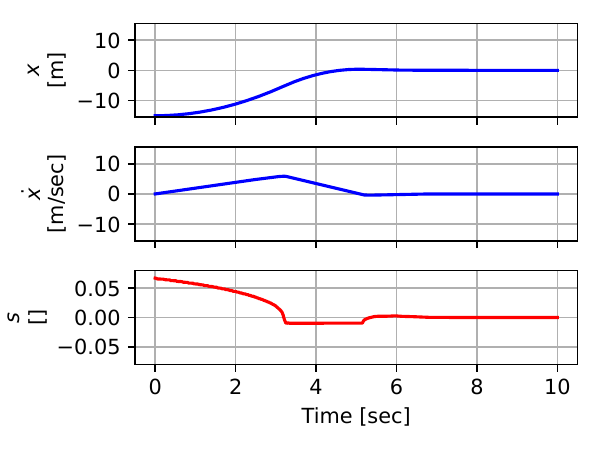}}
			
			\caption{Numerical results for context $c_f$.}
			\label{fig:c6c}
		\end{figure*}

\begin{figure*}[htp]
				\centering
				\vspace{-10pt}
				\subfloat[Feedback law $f$]{\includegraphics[width=6cm]{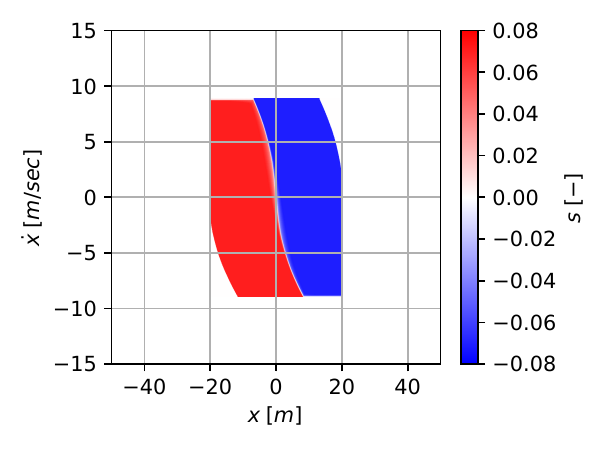}}
				\subfloat[Dimensionless feedback law $f^*$]{\includegraphics[width=6cm]{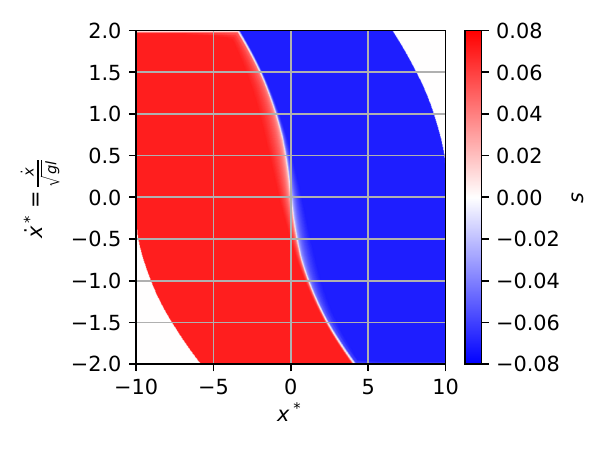}}   
				\subfloat[Optimal trajectory, starting at $x = -5 l$]{\includegraphics[width=6cm]{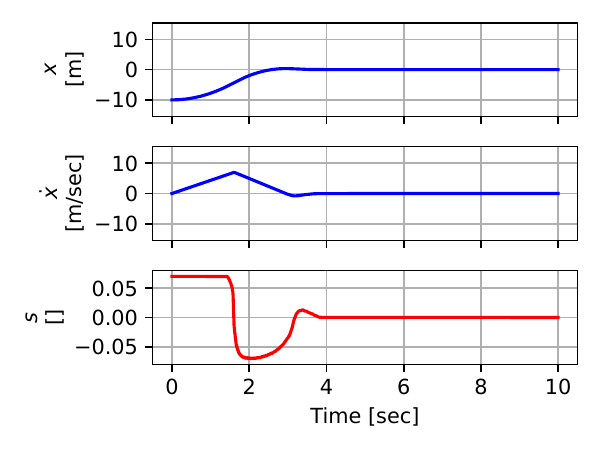}}
			\caption{Numerical results for context $c_g$.}
			\label{fig:c7c}
		\end{figure*}

\begin{figure*}[htp]
			
				\centering
				\vspace{-10pt}
				\subfloat[Feedback law $f$]{\includegraphics[width=6cm]{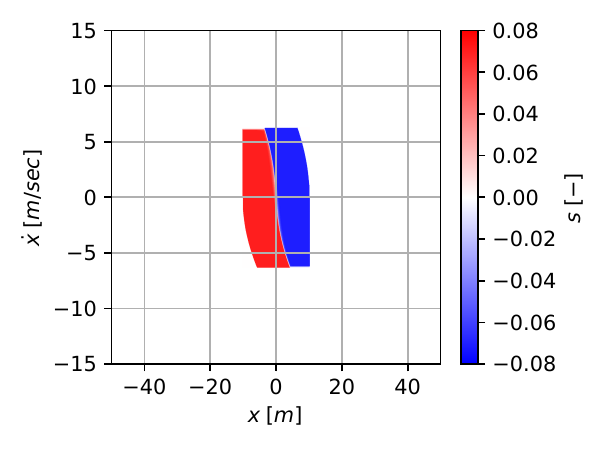}}
				\subfloat[Dimensionless feedback law $f^*$]{\includegraphics[width=6cm]{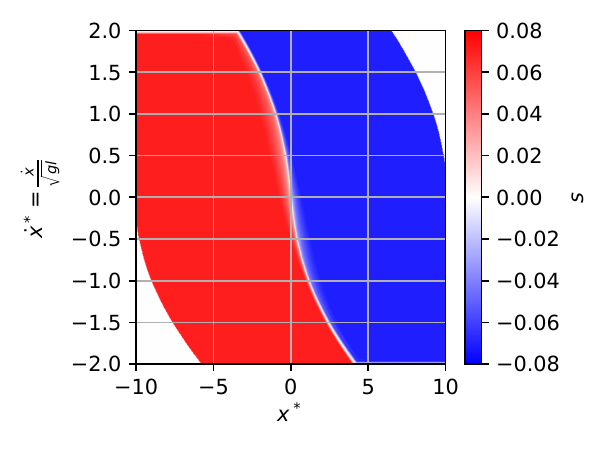}}   
				\subfloat[Optimal trajectory, starting at $x = -5 l$]{\includegraphics[width=6cm]{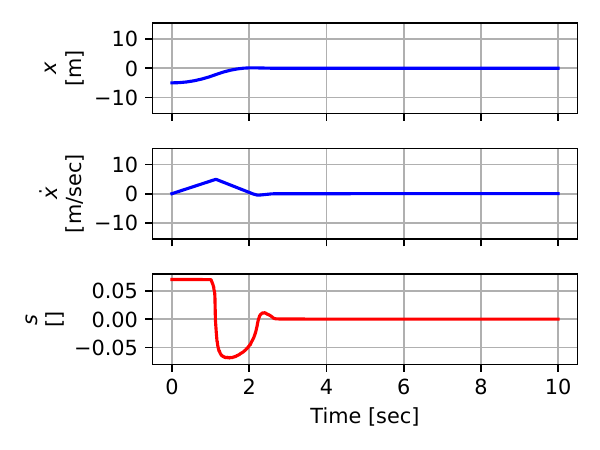}}
			
			\caption{Numerical results for context $c_h$.}
			\label{fig:c8c}
		\end{figure*}

\begin{figure*}[tp]
			
				\centering
				\vspace{-10pt}
				\subfloat[Feedback law $f$]{\includegraphics[width=6cm]{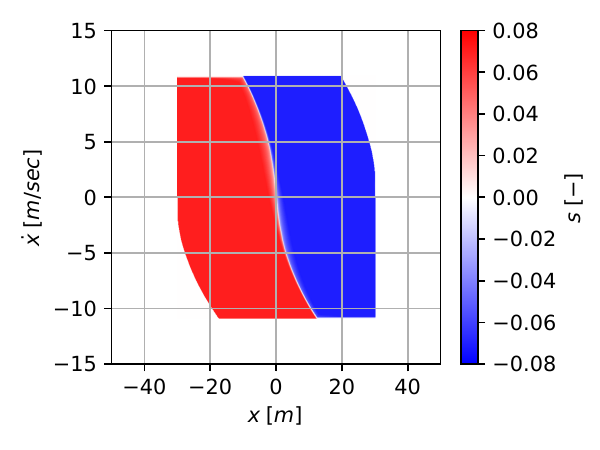}}
				\subfloat[Dimensionless feedback law $f^*$]{\includegraphics[width=6cm]{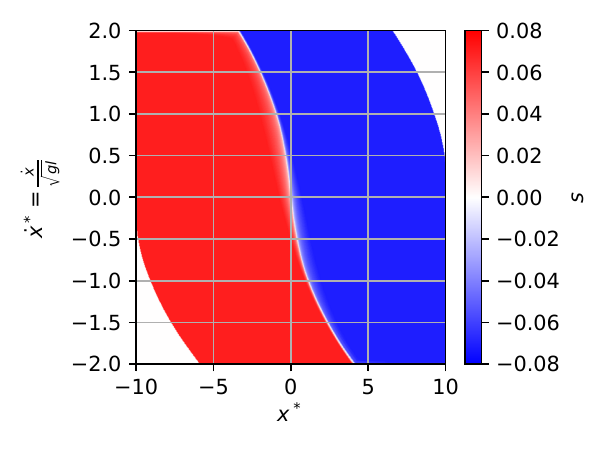}}   
				\subfloat[Optimal trajectory, starting at $x = -5 l$]{\includegraphics[width=6cm]{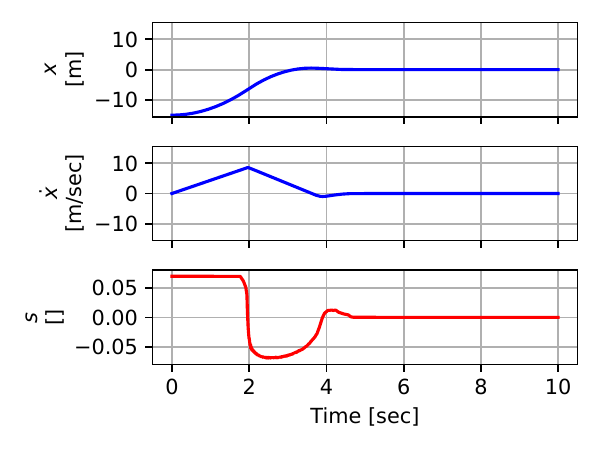}}
			
			\caption{Numerical results for context $c_i$.}
			\label{fig:c9c}
		\end{figure*}

Figures \ref{fig:c1c} to \ref{fig:c9c} illustrate that, for each subset with an equal dimensionless context, solutions are equal within each dimensionally similar subset when scaled into the dimensionless form. This was, again, the expected result predicted by the dimensional analysis presented in Section \ref{sec:dimenanalysis}. In terms of how to use this in a practical scenario, this exemplifies how various cars (which are different but which share the same ratios) could share a braking policy, for instance.
%%%%%%%%%%%%%%%%%%%%%%%%%%%%
		%%%%%%%%%%%%%%%%%%%%%%%%%%%%
		%%%%%%%%%%%%%%%%%%%%%%%%%%%%
		%%%%%%%%%%%%%%%%%%%%%%%%%%%%
		%%%%%%%%%%%%%%%%%%%%%%%%%%%%
		%%%%%%%%%%%%%%%%%%%%%%%%%%%
		%%%%%%%%%%%%%%%%%%%%%%%%%%%%
		%%%%%%%%%%%%%%%%%%%%%%%%%%%%
		%%%%%%%%%%%%%%%%%%%%%%%%%%%%
		%%%%%%%%%%%%%%%%%%%%%%%%%%%%

\subsubsection{Methodology}
\label{sec:metho2}

The same methodology as the pendulum example (see Section \ref{sec:metho}) was used for the car motion control problem. The time step was set to 0.025 s, the state space was discretized into an even 501 x 501 grid, and the continuous slip input was discretized into 101 discrete control options. Additional domain parameters were set as follows:
%%%%%%%%%%%%%%%%%%%%%%
\begin{align}
x_{max} &= 10l \\
\dot{x}_{max} &= \frac{ 2 }{\sqrt{gl}} \\
t_{max} &= 10 \frac{ x_{max} }{ v_{max} } 
\end{align}
%%%%%%%%%%%%%%%%%%%%%%
The source code is available online at the following link: \url{https://github.com/alx87grd/DimensionlessPolicies}, and this Google Colab page allows users to reproduce the results: \url{https://colab.research.google.com/drive/1-CSiLKiNLqq9JC3EFLqjR1fRdICI7e7M?usp=share_link}.

%%%%%%%%%%%%%%%%%%%%%%%%%%%%%%%%%%%%%%%%%%%%%%%%%
\section{Case studies with closed-form parametric policies}
\label{sec:closedfrom}

To better understand the concept of a dimensionless policy, in this section two examples based on well-known closed-form solutions to classical motion control problems are presented to illustrate how using Theorem \ref{theo:ab} can be equivalent to substituting new system parameters in an analytical solution.

%%%%%%%%%%%%%%%%%%%%%%%%%%%%%%%%%%%%%%%%%%%%%%%%%%%%%%%
\subsection{ Dimensionless linear quadratic regulator}
\label{sec:lqr}

The first example is based on the linear quadratic regulator (LQR) solution \cite{kalman_contributions_1960} for the linearized pendulum that allows for a closed-form analytical solution of optimal policy. This allow us to compared the method of transferring the policy with the proposed scaling law of Equation \eqref{eq:ab_transform}, to the method of transferring the policy by substituting the new system parameters in the analytical solution.

Here, we consider a simplified version of the pendulum swing-up problem (see Section \ref{sec:optimalswingup}) and a linearized version of the equation of motion is used, as follows:
%%%%%%%%%%%%%%%%%%%%%%
\begin{equation}
ml^2 \ddot{\theta} - mgl \theta = \tau
\label{eq:pendulum_dynamics_lqr}
\end{equation}
%%%%%%%%%%%%%%%%%%%%%%
The same infinite horizon quadratic cost function is used, as follows:
%%%%%%%%%%%%%%%%%%%%%%
\begin{equation}
J = \int_0^{\infty}{ \left( q^2 \theta^2 + \, \tau^2 \right) dt }
\label{eq:pendulum_cost_lqr}
\end{equation}
%%%%%%%%%%%%%%%%%%%%%%
However, no constraints on the torque are included in this problem. All parameters are also assumed to be time-independent constant. The same variables are used in this problem definition as before, except that the torque limit $\tau_{max}$ variable is absent. The global policy solution should then have the following form:
%%%%%%%%%%%%%%%%%%%%%%
\begin{equation}
\underbrace{\tau}_{\text{inputs}}
=
\pi_{lqr} \left(
\underbrace{ \theta, \dot{\theta} }_{\text{states}},
\underbrace{
\underbrace{ m , g , l }_{\text{system parameters}},
\underbrace{ q }_{\text{task parameters}}
}_{\text{context $c$}}
\right)
\label{eq:lqr_policy_form}
\end{equation}
We can thus select the same dimensionless $\Pi$ groups as in Section \ref{sec:dimanalpendulum} and conclude that Equation \eqref{eq:lqr_policy_form} can be restated under the following dimensionless form:
%%%%%%%%%%%%%%%%%%%%%%
\begin{equation}
\tau^*
=
\pi^*_{lqr}  \left(
 \theta^*, \dot{\theta}^* ,
 q^* 
\right)
\label{eq:lqr_dimpolicy_form}
\end{equation}

%%%%%%%%%%%%%%%%%%%%%%
\begin{Proposition}
For this motion control problem, defined by Equation \eqref{eq:pendulum_dynamics_lqr} and Equation \eqref{eq:pendulum_cost_lqr}, an analytical solution exists and the optimal policy is given by:
%%%%%%%%%%%%%%%%%%%%%%
\begin{align}
\tau &= 
\Biggl[ 
mgl +
\sqrt{ (mgl)^2 + q^2} \Biggr] \theta
\nonumber \\ &+
\Biggl[ \sqrt{ 2 ml^2 \Bigl( mgl+ \sqrt{ (mgl)^2 + q^2}} \Bigr) \Biggr] \dot{\theta}
\label{eq:lqr_policy}
\end{align}
%%%%%%%%%%%%%%%%%%%%%%
\end{Proposition}
\begin{proof}
See Appendix.% \ref{sec:lqr_proof}.
\end{proof}

Applying Equation \eqref{eq:f2star} to this feedback law leads to the dimensionless form, using $G=mgl$ and $H=ml^2$ for shortness, as follows:
%%%%%%%%%%%%%%%%%%%%%%
\begin{align}
\tau^* &= f^* ( x^* ) = \left[T_u(c)\right] f \left( \left[T_x^{-1}(c)\right] \; x^* \right) = \frac{1}{G} f( \theta^* ,\omega\dot{\theta}^* )  \\
\tau^* &= 
\left[\frac{1}{G}\right]  \left[ G + \sqrt{G^2+q^2}\right] \theta^*  
\nonumber \\ &+
\left[\frac{1}{G}\right] \left[\sqrt{2H \left( G+\sqrt{G^2+q^2} \right) }
\right] \left[\omega \dot{\theta}^*\right] \\
\tau^* &= 
\left[ 1 + \sqrt{\frac{G^2+q^2}{G^2}}\right] \theta^* 
\nonumber \\ &+
\left[\sqrt{\frac{2H\omega^2}{G}\frac{G+\sqrt{G^2+q^2}}{G}}
\right] \dot{\theta}^* \\
\tau^* &= 
\left[
1 + \sqrt{ 1 + (q^*)^2}
\right] \theta^*
+
\left[
\sqrt{2} \sqrt{ 1 + \sqrt{ 1 + (q^*)^2}}
\right] \dot{\theta}^*
\label{eq:dimlesslqr}
\end{align}
%%%%%%%%%%%%%%%%%%%%%%
The dimensionless policy is only a function of the dimensionless states and the dimensionless cost parameter $q^*$, as predicted by Equation \eqref{eq:lqr_dimpolicy_form} based on the dimensional analysis. It is interesting to note that Equation \eqref{eq:dimlesslqr} represents the core generic solution to the LQR problem and is independent of unit and scale.

We can also use this analytical policy solution to demonstrate Theorem \ref{theo:ab}, i.e. show that scaling the policy with Equation \eqref{eq:ab_transform} is equivalent to substituting new context variables when the contexts are dimensionally similar. 
%%%%%%%%%%%%%%%%%%%%%%
\begin{Proposition}
Suppose that we have two context instances, labeled $a$ and $b$, and that we use the global policy solution of Equation \eqref{eq:lqr_policy} to obtain two versions of context-specific feedback laws:
%%%%%%%%%%%%%%%%%%%%%%
\begin{align}
f_a(\theta,\dot{\theta}) &=
\left[ G_a + \sqrt{G_a^2+q_a^2}\right] \theta 
\nonumber \\ &+
\left[\sqrt{2H_a(G_a+\sqrt{G_a^2+q_a^2})}
\right] \dot{\theta} \label{eq:lqr_a}\\
f_b(\theta,\dot{\theta}) &=
\left[ G_b + \sqrt{G_b^2+q_b^2}\right] \theta 
\nonumber \\ &+
\left[\sqrt{2H_b(G_b+\sqrt{G_b^2+q_b^2})}
\right] \dot{\theta} \label{eq:lqr_b}
\end{align}
%%%%%%%%%%%%%%%%%%%%%%
where 
%%%%%%%%%%%%%%%%%%%%%%
\begin{align}
G_a=m_a g_a l_a \quad H_a=m_a l_a^2 \\
G_b=m_b g_b l_b \quad H_b=m_b l_b^2 
\end{align}
%%%%%%%%%%%%%%%%%%%%%%
Based on Theorem \ref{theo:ab}, if $q_a^* = q_b^*$ we can obtain $f_b$ directly by scaling $f_a$ based on Equation \eqref{eq:ab_transform} as follow:
%%%%%%%%%%%%%%%%%%%%%%
\begin{align}
f_b(\theta,\dot{\theta}) &= \left[ \frac{G_b}{G_a} \right] f_a \left( \theta , \left[ \frac{\omega_a}{\omega_b} \right] \dot{\theta} \right) \label{eq:lqr_ab}
\end{align}
%%%%%%%%%%%%%%%%%%%%%%
where 
%%%%%%%%%%%%%%%%%%%%%%
\begin{align}
\omega_a = \sqrt{G_a/H_a} \quad \omega_b = \sqrt{G_b/H_b}
\end{align}
%%%%%%%%%%%%%%%%%%%%%%
\end{Proposition}
\begin{proof}
If we substitute $f_a$ in Equation \eqref{eq:lqr_ab} by the analytical solution given by Equation \eqref{eq:lqr_a}, and then distribute the multiplying scaling factors we obtain:
%%%%%%%%%%%%%%%%%%%%%%
\begin{align}
f_b(\theta,\dot{\theta}) &= \left[ \frac{G_b}{G_a} \right] \Biggl( \left[ G_a + \sqrt{G_a^2+q_a^2}\right] \theta 
\nonumber \\ &+
\left[\sqrt{2H_a(G_a+\sqrt{G_a^2+q_a^2})}
\right]  \left[ \frac{\omega_a}{\omega_b} \right] \dot{\theta} \Biggr)
\\
f_b(\theta,\dot{\theta}) &= G_b \left[
1 + \sqrt{ 1 + (q_a^*)^2}
\right] \theta
\nonumber \\ &+
G_b
\left[
\sqrt{2} \sqrt{ 1 + \sqrt{ 1 + (q_a^*)^2}} \right] \frac{\dot{\theta}}{\omega_b} 
\\
% \end{align}
% %%%%%%%%%%%%%%%%%%%%%%
% where $q_a^* = q_a / G_a$ is the dimensionless version of the cost parameter. Then, distributing $G_b$ lead to:
% %%%%%%%%%%%%%%%%%%%%%%
% \begin{align}
% f_b(\theta,\dot{\theta}) &= \left[
% G_b + \sqrt{ G_b^2 + (G_b q_a^*)^2}
% \right] \theta
% +
% \left[
% \sqrt{\frac{2G_b}{\omega_b^2}} \sqrt{ G_b + \sqrt{ G_b^2 + (G_b q_a^*)^2}} \right] \dot{\theta}
% \\
f_b(\theta,\dot{\theta}) &= \left[
G_b + \sqrt{ G_b^2 + (G_b q_a^*)^2}
\right] \theta
\nonumber \\ &+
\left[
\sqrt{2 H_b} \sqrt{ G_b + \sqrt{ G_b^2 + (G_b q_a^*)^2}} \right] \dot{\theta}
\end{align}
%%%%%%%%%%%%%%%%%%%%%%
which is equivalent to Equation \eqref{eq:lqr_b} when
%%%%%%%%%%%%%%%%%%%%%%
\begin{align}
G_b q_a^* = q_b \quad \text{or equivalently } \quad q_a^* = q_b^*
\end{align}
%%%%%%%%%%%%%%%%%%%%%%
which is the condition of having equal dimensionless contexts ($c_a^* = c_b^*$) for this motion control problem. 
\end{proof}

This example illustrates that applying the scaling of Equation \eqref{eq:ab_transform} based on the dimensional analysis framework is equivalent to changing the context variables in an analytical solution when the dimensionless context variables are equal.

\subsection{Dimensionless computed torque}

The second example is again based on the pendulum, but using the computed torque control technique \cite{asada_robot_1986}. This also allow us to compared the method of transferring the policy with the proposed scaling law of Equation \eqref{eq:ab_transform}, to the method of transferring the policy by substituting the new system parameters in the analytical solution. This example is not based on a quadratic cost function, as opposed to previous examples, to illustrate the flexibility of the proposed schemes.

Here, we present a second analytical example, A computed torque feedback law is a model-based policy (assuming that there are no torque limits) that is the solution to the motion control problem of making a mechanical system that converges on a desired trajectory, with a specified second-order exponential time profile defined by the following equation:
%%%%%%%%%%%%%%%%%%%%%%
\begin{equation}
0 = (\ddot{\theta}_d - \ddot{\theta})+ 2 \omega_d \zeta (\dot{\theta}_d - \dot{\theta}) + \omega_d^2 (\theta - \theta)
\end{equation}
%%%%%%%%%%%%%%%%%%%%%%
For the specific case of the pendulum-swing up problem, we assume that all parameters are time-independent constants and that our desired trajectory is simply the upright position ($\ddot{\theta}_d = \dot{\theta}_d = \theta_d = 0$), leaving only two parameters to define the tasks: $\omega_d$ and $\zeta$. Then, the computed torque policy takes the following form:
%%%%%%%%%%%%%%%%%%%%%%
\begin{align}
\underbrace{\tau}_{\text{input}} 
&=
\pi_{ct} \left(
\underbrace{ \theta, \dot{\theta} }_{\text{states}},
\underbrace{
\underbrace{ m , g , l }_{\text{system parameters}},
\underbrace{ \omega_d , \zeta }_{\text{task parameters}}
}_{\text{context $c$}}
\right) 
%\label{tb:ct2}
\end{align}
%%%%%%%%%%%%%%%%%%%%%%
and the analytical solution is as follows:
%%%%%%%%%%%%%%%%%%%%%%
\begin{align}
\tau &= mgl \sin \theta - 2 m l^2 \omega_d \zeta \dot{\theta} - m l^2 \omega_d^2 \theta \label{eq:ct}
\end{align}
%%%%%%%%%%%%%%%%%%%%%%
Here, the context includes the system parameters and two variables characterizing the convergence speed. Note that the task parameters directly define the desired behavior, as opposed to the previous examples where they were defining the behavior indirectly thought a cost function. The states, control inputs, and system parameters are the same as before; only the task parameters differ, and their dimensions are presented in Table \ref{tb:ctc}.
%%%%%%%%%%%%%%%%%%%%%%%%%%%%%%%%%%%%%%%%%%%%
\begin{table}[htb]
    \centering % center the table
    \caption{Computed torque task variables.} 
    \label{tb:ctc}
    \begin{tabular}{p{1.0cm} p{2.5cm} p{1.0cm} p{1.5cm}}
        \hline \hline \noalign{\smallskip} \noalign{\smallskip} \noalign{\smallskip} \noalign{\smallskip}
        %%%%%%%%%%%%%%%%%%%%%%
        \textbf{Variable} & \textbf{Description} & \textbf{Units} & \textbf{Dimensions} \\ 
        %    %%%%%%%%%%%%%%%%%%%%%%
        %    \hline \hline \noalign{\smallskip} 
        %    \multicolumn{4}{c}{\textbf{Control inputs}}\\ \noalign{\smallskip}  \hline \hline
        %    \noalign{\smallskip} 
        %    %%%%%%%%%%%%%%%%%%%%%%
        %    $\tau$ & Actuator torque & $Nm$ & [$ML^2T^{-2}$]\\ 
        %    %%%%%%%%%%%%%%%%%%%%%%
        %    \hline \hline \noalign{\smallskip} 
        %    \multicolumn{4}{c}{\textbf{State variables}}\\ \noalign{\smallskip}  \hline \hline \noalign{\smallskip} 
        %    %%%%%%%%%%%%%%%%%%%%%%
        %    $\theta$ & Joint angle & $rad$ & []\\ \noalign{\smallskip} \hline \noalign{\smallskip}
        %    $\dot{\theta}$ & Joint angular velocity & $rad/sec$ & [$T^{-1}$] \\
        %    %%%%%%%%%%%%%%%%%%%%%%
        %    \hline \hline \noalign{\smallskip} 
        %    \multicolumn{4}{c}{\textbf{System parameters}}\\ \noalign{\smallskip}  \hline\hline  \noalign{\smallskip} 
        %    %%%%%%%%%%%%%%%%%%%%%%
        %    $mgl$ & Maximum gravitational torque  & $Nm$ & [$ML^2T^{-2}$]  \\ \noalign{\smallskip} \hline \noalign{\smallskip}
        %    $\omega = \sqrt{\frac{g}{l}}$ & Natural frequency & $sec^{-1}$ & [$T^{-1}$]  \\ \noalign{\smallskip} \hline \noalign{\smallskip}
        % %%%%%%%%%%%%%%%%%%%%%%
        \hline \hline \noalign{\smallskip} 
        \multicolumn{4}{c}{\textbf{Task parameters}}\\ \noalign{\smallskip}  \hline\hline  \noalign{\smallskip} 
        %%%%%%%%%%%%%%%%%%%%%%
        $\omega_d$ & Desired closed-loop frequency & $s^{-1}$ & [$T^{-1}$]  \\ \noalign{\smallskip} \hline \noalign{\smallskip}
        $\zeta$ &  Desired closed-loop damping & $-$ & [-]  \\ \noalign{\smallskip} \hline \noalign{\smallskip}
    \end{tabular}
\end{table}
%%%%%%%%%%%%%%%%%%%%%%%%%%%%%%%%%%%%%%%%%%%%

Here, seven variables and only $p=2$ independent dimensions ( $ML^2T^{-2}$ and $T^{-1}$ ) are involved. Thus, five dimensionless groups can be formed, as follows:
%%%%%%%%%%%%%%%%%%%%%%
\begin{equation}
1 + (n=2) + (m=4) - ( p = 2 ) = 5
\end{equation}
%%%%%%%%%%%%%%%%%%%%%%
Using $mgl$ and $\omega$, the system parameters, as the repeating variables leads to the following dimensionless groups:
%%%%%%%%%%%%%%%%%%%%%%
\begin{align}
\Pi_1 &= \tau^* = \frac{\tau}{mgl} \quad \quad \frac{[ML^2T^{-2}]}{[M][LT^{-2}][L]} \\
\Pi_2 &= \theta^* = \theta \quad \quad [-]\\
\Pi_3 &= \dot{\theta}^* = \frac{ \dot{\theta}  }{ \omega } \quad \quad \frac{[T^{-1}]}{[T^{-1}]} \\
\Pi_4 &= \omega_d^* = \frac{\omega_d}{\omega} \quad \quad \frac{[T^{-1}]}{[T^{-1}]} \\
\Pi_5 &= \zeta^* = \zeta \quad \quad [-]
\end{align}
%%%%%%%%%%%%%%%%%%%%%%
Then, applying the Buckingham $\pi$ theorem tells us that the computed torque policy can be restated as the following relationship between the dimensionless variables:
%%%%%%%%%%%%%%%%%%%%%%
\begin{equation}
\tau^*
=
\pi^*_{ct} \left(
\theta, \dot{\theta}^*,
\omega_d^* , \zeta^* 
\right)
\label{eq:ct_dim_rel}
\end{equation}
%%%%%%%%%%%%%%%%%%%%%%

Here, we can confirm directly (since we have an analytical solution) that applying Equation \eqref{eq:f2star} to the computed torque feedback law given by Equation \eqref{eq:ct} leads to the following dimensionless form:
%%%%%%%%%%%%%%%%%%%%%%
\begin{align}
\tau^* &= \left[ \frac{1}{mgl} \right] \left( mgl \sin \theta - 2 m l^2 \omega_d \zeta \left( \omega \dot{\theta}^*\right) - m l^2 \omega_d^2 \theta \right) \\
\tau^*
&=
\sin \theta^*
- 2 \omega_d^* \zeta^* \dot{\theta}^* 
- (\omega_d^*)^2 \theta^*
\end{align}
%%%%%%%%%%%%%%%%%%%%%%
thereby confirming the structure predicted by Equation \eqref{eq:ct_dim_rel} based on the dimensional analysis.

We can, again, use this example to demonstrate Theorem \ref{theo:ab} and show that, when the dimensionless context is equal, scaling a policy using Equation \eqref{eq:ab_transform} is equivalent to substituting new values of the system parameters into the analytical equation. 
%%%%%%%%%%%%%%%%%%%%%%
\begin{Proposition}
Suppose that we have two context instances, labeled $a$ and $b$, and that we use the global policy solution of Equation \eqref{eq:ct} to obtain two versions of context-specific feedback laws:
%%%%%%%%%%%%%%%%%%%%%%
\begin{align}
f_a(\theta,\dot{\theta}) &= G_a \sin \theta - 2 H_a \omega_{d,a} \zeta_a \dot{\theta} - H_a \omega_{d,a}^2 \theta \label{eq:ctc_a} \\
f_b(\theta,\dot{\theta}) &= G_b \sin \theta - 2 H_b \omega_{d,a} \zeta_b \dot{\theta} - H_b \omega_{d,a}^2 \theta \label{eq:ctc_b}
\end{align}
%%%%%%%%%%%%%%%%%%%%%%
Based on Theorem \ref{theo:ab}, if $\omega^*_{d,a}=\omega^*_{d,b}$ and $\zeta_a^* = \zeta_b^*$ we can obtain $f_b$ directly by scaling $f_a$ based on Equation \eqref{eq:ab_transform} as follow:
%%%%%%%%%%%%%%%%%%%%%%
\begin{align}
f_b(\theta,\dot{\theta}) &= \left[ \frac{G_b}{G_a} \right] f_a \left( \theta , \left[ \frac{\omega_a}{\omega_b} \right] \dot{\theta} \right) \label{eq:ctc_ab}
\end{align}
%%%%%%%%%%%%%%%%%%%%%%
\end{Proposition}
\begin{proof}
If we substitute $f_a$ in Equation \eqref{eq:ctc_ab} by the analytical solution given by Equation \eqref{eq:ctc_a}, and then distribute the multiplying scaling factors we obtain:
%%%%%%%%%%%%%%%%%%%%%%
\begin{align}
f_b(\theta,\dot{\theta}) &= \left[ \frac{G_b}{G_a} \right] \Bigr[  G_a \sin \theta - 2 H_a \omega_{d,a} \zeta_a \left[\frac{\omega_a}{\omega_b} \right] \dot{\theta}  
\nonumber \\ &-
H_a \omega_{d,a}^2 \theta \Bigl]  
\\
f_b(\theta,\dot{\theta}) &= G_b \left[  \sin \theta - 2 \frac{\omega_{d,a} }{\omega_a } \zeta_a \frac{\dot{ \theta }}{\omega_b} - \left(\frac{ \omega_{d,a} }{\omega_a }\right)^2 \theta \right]
\\
f_b(\theta,\dot{\theta}) &= G_b  \sin \theta - 2 H_b  
%\underbrace{
\left( \frac{\omega_b}{\omega_a} \omega_{d,a} 
\right) \zeta_a \dot{ \theta }
%}_{(\omega_d \zeta )_b}
\nonumber \\ &-
H_b 
%\underbrace{
\left( 
\frac{\omega_b}{\omega_a} \omega_{d,a}
\right)^2 
%}_{(\omega_d)_b^2 }
\theta
\end{align}
%%%%%%%%%%%%%%%%%%%%%%
which is exactly equivalent to Equation \eqref{eq:ctc_b} (i.e., equivalent to substituting the $a$ instance of the context variables to the $b$ instance) if:
%%%%%%%%%%%%%%%%%%%%%%
\begin{align}
\frac{\omega_b}{\omega_a} \omega_{d,a} = \omega_{d,b} 
\quad \text{and}  \quad \zeta_a = \zeta_b
\end{align}
%%%%%%%%%%%%%%%%%%%%%%
which is the dimensional similarity condition ($c_a^* = c_b^*$) for this motion control problem:
%%%%%%%%%%%%%%%%%%%%%%
\begin{align}
\frac{\omega_{d,a}}{\omega_a}  = \omega^*_a =  \omega^*_b = \frac{\omega_{d,b}}{\omega_b} 
\quad \text{and}  \quad \zeta_a^* = \zeta_b^*
\end{align}
%%%%%%%%%%%%%%%%%%%%%%
\end{proof}

%%%%%%%%%%%%%%%%%%%%%%%%%%%%%%%%%%%%%%%%%%%%%%%%%%%%
\section{Conclusion}

The dimensional analysis of physically meaningful control policies, leveraging the Buckingham $\pi$ theorem, leads to two interesting theoretical results: \textbf{1) }In dimensionless form, the solution to a motion control problem involves a reduced number of parameters. \textbf{2)} It is possible to exactly transfer a feedback law between similar systems without any approximation, simply by scaling the input and output of any type of control law appropriately, including via numerically generated black box mapping. However, the main practical limitation of this approach is that if the condition of dimensional similarity ($c_a^*=c_b^*$) is not met exactly, then there is no theoretical guarantees regarding whether a policy is transferable without additional assumptions, as the discussed concept of regimes of behaviour. Also, we demonstrated how those results can be used to transfer exactly even discontinuous black-box policies between similar systems, using two simple examples of dynamical systems and numerically generated optimal feedback laws. An interesting direction for further exploration would be investigating how good an approximation is when a feedback law is transferred from a context that is not exactly similar but close. Also, it would be interesting to test the concept of dimensionless policies to empower a reinforcement learning scheme that could collect data from various, but dimensionally similar, systems to accelerate the learning process.

\appendix
%\section*{LQR Analytic Solution}
%\label{sec:lqr_proof}
%%%%%%%%%%%%%%%%%%%%%%%%%%%%%%%%%%%%%%%%%%%%%%%%%%%%%%%%%%%%%%%%%%
In this section, we show that the policy given by Equation \eqref{eq:lqr_policy} is optimal with respect to the LQR problem defined in Section \ref{sec:lqr}. We can write the equation of motion given by Equation \eqref{eq:pendulum_dynamics_lqr} in state-space form, using $G=mgl$ and $H=ml^2$, as follows:
%%%%%%%%%%%%%%%%%%%%%%
\begin{equation}
\frac{d}{dt} 
\begin{bmatrix}
    \theta \\
    \dot{\theta} 
\end{bmatrix}
= 
\underbrace{
\begin{bmatrix}
    0 & 1 \\
    G/H & 0 
\end{bmatrix}
}_{A}
\underbrace{
\begin{bmatrix}
    \theta \\
    \dot{\theta} 
\end{bmatrix}
}_{x}
+
\underbrace{
\begin{bmatrix}
    0 \\
    1/H
\end{bmatrix}
}_{B}
\underbrace{
\begin{bmatrix}
    \tau
\end{bmatrix}
}_{u}
\end{equation}
%%%%%%%%%%%%%%%%%%%%%%

Then, by adapting a solution from \cite{hanks_closed-form_1991}, if we parameterize the weight matrix of the cost function as follows:
%%%%%%%%%%%%%%%%%%%%%%
\begin{equation}
J = \int_0^{\infty}{ \left( x^T 
\underbrace{
\begin{bmatrix}
    a(a-2G) & 0 \\
    0       & b^2 - 2 a H
\end{bmatrix}
}_{Q} x
+
u^T 
\underbrace{
\begin{bmatrix}
    1 
\end{bmatrix}
}_{R} u
\right) dt } 
\end{equation}
%%%%%%%%%%%%%%%%%%%%%%
the optimal cost-to-go is given by:
%%%%%%%%%%%%%%%%%%%%%%
\begin{equation}
J =   x^T 
\underbrace{
\begin{bmatrix}
    b(a-G) & aH \\
    aH       & bH
\end{bmatrix}
}_{S} x
\end{equation}
%%%%%%%%%%%%%%%%%%%%%%
and the optimal feedback policy is given by:
%%%%%%%%%%%%%%%%%%%%%%
\begin{equation}
u = 
- \underbrace{\left[ R^{-1} B^T S\right]}_{K} x
=
-
\underbrace{
\begin{bmatrix}
    a & b \\
\end{bmatrix}
}_{K} x
\label{eq:kx}
\end{equation}
%%%%%%%%%%%%%%%%%%%%%%
This solution can by verified by substituting matrices into the algebraic Riccati equation given by:
%%%%%%%%%%%%%%%%%%%%%%
\begin{equation}
0 = SA + A^T S - SBR^{-1}B^TS + Q
\end{equation}
%%%%%%%%%%%%%%%%%%%%%%
since the problem fits into the framework of the classical infinite horizon LQR result \cite{bertsekas_dynamic_2012}. Then, we can see that the cost function defined in Section \ref{sec:lqr} is a special case, where $Q_{11} = q^2$ and $Q_{22}=0$, leading to the following equations:
%%%%%%%%%%%%%%%%%%%%%%
\begin{align}
q^2 &= a(a-2G) \\
0   &= b^2 -2aH
\end{align}
%%%%%%%%%%%%%%%%%%%%%%
Solving for $a$ and $b$, and retaining the positive solution, leads to the following:
%%%%%%%%%%%%%%%%%%%%%%
\begin{align}
a &= G + \sqrt{G^2+q^2} \\
b &= \sqrt{2aH} = \sqrt{2H \left(G + \sqrt{G^2+q^2} \right)}
\end{align}
%%%%%%%%%%%%%%%%%%%%%%
which, when substituted into Equation \eqref{eq:kx}, is equal to the policy given by Equation \eqref{eq:lqr_policy} in Section \ref{sec:lqr}.

%\newpage

\bibliographystyle{IEEEtran}
\bibliography{zotero_alex}

% that's all folks
\end{document}